\documentclass[11 pt]{article}

\usepackage[english]{babel}

\usepackage[letterpaper,top=2cm,bottom=2cm,left=3cm,right=3cm,marginparwidth=1.75cm]{geometry}

\usepackage[title]{appendix}

\usepackage{amsmath}
\usepackage{graphicx}
\usepackage[colorlinks=true, allcolors=blue]{hyperref}

\usepackage{pdflscape} 
\usepackage{afterpage}

\usepackage[utf8]{inputenc}
\usepackage{amsmath, amssymb,mathtools}
\usepackage{algorithm}
\usepackage{algpseudocode}
\usepackage{float}
\usepackage[caption = false]{subfig}
\usepackage{graphicx,appendix}
\usepackage{bm}
\usepackage{tikz}

\usetikzlibrary{shapes,arrows}

\usepackage{lineno}
\usepackage{wrapfig}

\usepackage{colortbl}
\usepackage{xcolor}
\usepackage{soul}
\usepackage{comment}

\usepackage{amsthm}

\newtheorem{remark}{Remark}

\newcommand{\R}{\mathbb{R}}
\newcommand{\M}{\mathbb{M}}

\newcommand{\X}{\mathcal{X}}

\newcommand{\T}{\mathcal{T}}

\newcommand{\tr}{\text{tr}}
\newcommand{\cmt}[1]{} 
\newcommand{\bu}{\mathbf{u}}
\newcommand{\by}{\mathbf{y}}
\newcommand{\m}{\mathbf{m}}
\newcommand{\bx}{\mathbf{x}}
\newcommand{\br}{\mathbf{r}}
\newcommand{\bU}{\mathbf{U}}
\newcommand{\F}{\mathbf{F}}
\newcommand{\bP}{\mathbf{P}}
\newcommand{\bX}{\mathbf{X}}
\newcommand{\bY}{\mathbf{Y}}
\newcommand{\bK}{\mathbf{K}}
\newcommand{\bH}{\mathbf{H}}
\newcommand{\bQ}{\mathbf{Q}}
\newcommand{\bR}{\mathbf{R}}
\newcommand{\bS}{\mathbf{S}}
\newcommand{\bI}{\mathbf{I}}
\newcommand{\be}{\mathbf{e}}

\newcommand{\Rc}[0]{\mathbf{R}} 
\newcommand{\Sc}[0]{\mathbf{Q}} 

\title{Robust Adaptive Meshing, Mesh Density Functions, and \\
Nonlocal Observations for Ensemble Based Data Assimilation%
\thanks{This work was supported in part by Office of Naval Research Grant N000142412218.}
}
\author{Jeremiah Buenger\thanks{Department of Mathematics, University of Kansas, Lawrence, Kansas, USA (jeremiah.buenger@ku.edu).},
\and Weizhang Huang\thanks{Department of Mathematics, University of Kansas, Lawrence, Kansas, USA (whuang@ku.edu).},
\and Erik S. Van Vleck\thanks{Department of Mathematics, University of Kansas, Lawrence, Kansas, USA (erikvv@ku.edu).}}

\begin{document}
\maketitle

\begin{abstract}
Adaptive spatial meshing has proven invaluable for the accurate, efficient computation of solutions of time dependent partial differential equations. In a DA context the use of adaptive spatial meshes addresses several factors that place increased demands on meshing; these include the location and relative importance of observations and the use of ensemble solutions. To increase the efficiency of adaptive meshes for data assimilation, robust look ahead meshes are developed that fix the same adaptive mesh for all ensemble members for the entire time interval of the forecasts and that incorporates the observations at the next analysis time. This allows for increased vectorization of the ensemble forecasts while minimizing interpolation of solutions between different meshes. The techniques to determine these robust meshes are based upon combining metric tensors or mesh density functions to define nonuniform meshes. We illustrate the robust ensemble look ahead meshes using traveling wave solutions of a bistable reaction-diffusion equation. Observation operators based on convolution type integrals and their associated metric tensors are derived. These further the goals of  
making efficient use of adaptive meshes in ensemble based DA techniques, developing and employing robust meshes that are effective for a range of similar behaviors in both the ensembles and the observations, and the integration with advanced numerical PDE techniques (a quasi-Lagrangian moving mesh DG technique employing embedded pairs for time stepping). Numerical experiments with different observation scenarios are presented for  a 2D inviscid Burgers' equation, a multi-component system, a 2D Shallow Water model, and for a coupled system of two 1D Kuramoto-Sivashinsky equations.

\end{abstract}

\section{Introduction}





Several challenges as well as opportunities exist for combining data assimilation (DA) with adaptive meshing techniques; techniques that have proven to be quite successful for the numerical solution of partial differential equations. In a DA context, other sources beyond accurate, efficient solution of forecast ensembles may require additional, variable resolution. Primary among this are the location and relative significance of observational data. Other factors such as the need to focus on important events and locations motivate the need for variable resolution in different space/time regions and for specific model variables. 
While natural that the physical model solutions and observational data are primary drivers for meshing, combining the requirements for both may enhance the skill of the DA. 
Even with the increased use of surrogate physical and data models, the need for flexible means of increasing resolution that are both accurate and efficient persists.
Detecting finer details is essential to understanding important real-time behavior.

Our contributions in this paper are to combine advanced techniques for numerical solution of partial differential equations with the development of algorithms for ensemble based DA. These algorithms are designed to take advantage of adaptive meshing techniques. 
Within each forecast/analysis DA cycle are two natural places to update mesh(es). One is before the analysis with an eye toward optimizing the analysis. Often this is done by determining a mesh common to all ensembles that also supports the observational data. In this way the analysis update is performed much like with an ordinary differential equation or a partial differential equation (PDE) on a fixed mesh, see, e.g., \cite{KHMVVZ22}. In this paper we consider another natural place to remesh, after the analysis and before the next set of ensemble forecast. 
We develop a common forecast mesh that also takes into account/optimize the subsequent analysis step. The construction of these meshes is based on metric tensors or mesh density functions. 
In particular, we consider non-stationary, nonlocal observation operators based on convolution integrals and derive a metric tensor to increase resolution in a neighborhood of the support of these observation operators. In addition, we employ a quasi-Lagrangian variational mesh adaptation technique that employs the derived metric tensors, combined with discontinuous Galerkin (DG) discretizations and
the use of strongly stability preserving (SSP) Runge-Kutta time stepping techniques with embedded pairs to facilitate adaptive time stepping.

There are several recent works on integrating adaptive spatial meshing techniques with DA, although most of the focus has been on PDE models in one space dimension (1D). This includes methods based on evolving meshes based on the solution of a differential equation, methods in which meshes are updated statically based upon interpolation, and remeshing techniques that add or subtract mesh points as the solution structure changes.
%
One approach is to interpolate the ensemble solutions to a common mesh at each observational time step and assimilate the PDE variables only. The common mesh approach is used here, as well as in \cite{Colin,Du,KHMVVZ22}. Another method is to assimilate both the PDE variables and the common mesh locations, as done in {\color{blue} \cite{Bonan, sampson2020}}.
In \cite{Bonan,sampson2020}, the state variables of the PDE were augmented with the position of the nodes and incorporated into a DA scheme. The test problem consisted of a two-dimensional ice sheet assumed to be radially symmetric; therefore, it reduced to a problem with one spatial dimension.
In \cite{Colin,Du} common meshes were developed based on combining the ensemble meshes through interpolation. This allowed for updating the mean and covariance for Kalman filter based DA techniques while allowing each ensemble member to evolve on its own independent mesh. That is, at each observational timestep, the ensemble members were interpolated to the common mesh, updated with the DA analysis, and then interpolated back to their respective meshes. Specifically, a uniform, non-conservative (remeshing allowed with the number of mesh points potentially varying with time) mesh was used in \cite{Colin}, with Lagrangian observations in one spatial dimension. Higher spatial dimensions were used in \cite{Du}, with a fixed common mesh refined near observation locations. In \cite{KHMVVZ22} a metric tensor or mesh density function characterization of nonuniform meshes was employed. This was then used to form meshes based on the combination of mesh density functions, one based on the ensemble forecasts and another based on the observation locations.
This provides increased flexibility that allows for the inclusion of different factors that impact the grid resolution, in particular by combining the mesh for ensemble forecasts with a mesh that focuses on the location of observations. The work \cite{sampson2020} employs a 1D non-conservative adaptive meshing scheme as in \cite{Colin} and extends this approach through the use of an adaptive common mesh, where, like in \cite{Bonan}, the state vector is augmented with the node locations.



This paper is outlined as follows. In Section \ref{background} we present the background on data assimilation, adaptive meshing, and outline the robust meshing strategy developed in this paper based on combining meshes through their metric tensors or mesh density functions. Section \ref{RobustMeshes} details the use of metric tensor intersection to create robust meshes, as well as strategies for obtaining robust look ahead meshes. Common meshes for forecasts and analysis are developed and robust look ahead meshes are illustrated using traveling wave solutions of a bistable reaction-diffusion equation. In Section \ref{EnsembleDA} ensemble based DA techniques are considered that potentially benefit from our developments and metric tensors/mesh density functions are employed for domain localization. 
Observation operators based upon a type of convolution integral are developed and solution dependent metric tensors are derived to concentrate mesh points for this type of observation operator. Section \ref{DACompParams} contains details of the computational techniques and data assimilation parameters to be used in our numerical experiments. We consider two hyperbolic PDEs, an inviscid Burgers' equation and a three component shallow water equation, both in two space dimensions. We employ DG methods and SSP Runge-Kutta time stepping techniques with adaptive time stepping using embedded pairs. In Section \ref{Numerics} we present numerical experiments for both hyperbolic equations as well as for a system of two coupled 1D Kuramoto-Sivashinsky equations. The conclusion in Section~\ref{Conclusion} contains some observations on the current work and directions for future research. 




\section{Background}\label{background}

Data assimilation techniques combine model forecasts and data to improve predictions and to quantify uncertainties typically in a Bayesian context (see, e.g., \cite{Carrassi2018,LSZ,ReichCotter15,VanLeeuwen18}). We consider a finite dimensional discrete time system for forecasts (a time dependent PDE discretized in space and time), for observation times
for $n=0,1,...,$
\begin{equation}\label{DAsetup}
    \left\{
    \begin{array}{rl}
    & \bu_{n+1} = \Psi(\bu_n) + \xi_n ,\cr
    & \by_{n+1} = {\mathcal H}(\bu_{n+1}) + \eta_{n},
    \end{array}
    \right.
\end{equation}
where, e.g., $\xi_n\sim N(0,\bQ),$ $\eta_{n}\sim N(0,\bR),$ 
and $u_0\sim N(\bu_0^b,\bP_0^b)$. Here ${\mathcal H}$ denotes the observation operator and $\by_{n+1}$ is the observational data at time $n+1.$ We denote by $N_y$ the dimension of the observation space so that $\by_{n+1}\in\R^{N_y}$, $N_u$ the dimension of the model space so $\bu_{n+1}\in\R^{N_u}$, and $N_e$ denotes the number of ensemble members. We are 
motivated by the typical case in which $N_e \ll N_y \ll N_u.$


Using a single forecast-analysis cycle to illustrate, meshes are potentially updated both for forecasts and for analysis. Here we denote by 
$\bU^f_n$ and $\bU^a_n$ the forecast and analysis ensembles, respectively, at time $n$. Referring to (\ref{eq:DApred}) and (\ref{eq:DAanalysis}) below
and \cite{KHMVVZ22}, the A options were employed and each forecast ensemble updates its own dedicated mesh. Before the analysis the ensemble meshes are combined with a mesh that focuses on the location of observations to form a common mesh that supports the analysis. The combining of meshes is accomplished by combining the metric tensors or mesh density functions that correspond to the different meshes. In this work, the B options are employed and a single robust, nonuniform common mesh is employed for all ensemble forecasts and for the analysis update. This is accomplished by using a small ensemble of forecasts that provide information on meshing requirements for the ensemble forecasts between observation times as well as the observations and their locations at the time of the next analysis update. In this way there is a single common mesh for all forecast ensemble members, ideally over the entire forecast window, and the common mesh is used for the next analysis update.

\begin{align}
    &\text{(Forecast)} \quad \begin{cases}
        \quad 
        &\text{{\bf A. Adaptive Forecast Meshes for each Ensemble Member;}} \\
        &\text{{\bf B. Small Ensemble to create Common Look Ahead Mesh, }}\\
        &\text{{\bf $\,\,\,$\quad Full Ensemble Forecast using the Look Ahead Mesh;}}
    \\
    &\,\qquad\qquad\qquad\qquad \bU_{n+1}^{f} = \F(\bU_n^{a})+\xi_n,
    \\
    \end{cases}\label{eq:DApred}\\
    &\text{(Analysis)} \quad \begin{cases}
    \quad &\text{{\bf A. Analysis Mesh Common to all Ensemble Members;}}\\
    \quad &\text{{\bf B. Forecast Mesh designed to Support Analysis;}}
    \\
    &\,\qquad\qquad\bU_{n+1}^a = 
    \bU_{n+1}^f + \mathbf{K}_{n+1}(\by_{n+1} - \mathcal{H} (\bU_{n+1}^f)),
    \end{cases}\label{eq:DAanalysis}
\end{align}

The goals for the mesh(es) supporting the forecasts and the analysis are different but intertwined. It is natural to remesh after the analysis step (e.g., when there is large change in $\bU$) and before the forecast step (to optimize accuracy, efficiency, etc., of both the forecasts and the DA skill). There is good potential to use the optimization problem the analysis is based upon to more directly develop improved mesh(es) to support the analysis.

When remeshing before the analysis, one has full knowledge of ensemble forecasts and their corresponding meshes. By combining the ensemble meshes with a mesh focused on the observations and any other important factors, a common mesh to all ensembles may be formed that supports the analysis. When remeshing after the current analysis and before the subsequent forecasts, a mesh update is often necessary since there may be large changes in the ensemble solutions as a result of the analysis. By forming a common mesh for the ensemble forecasts, there is a possibility to improve the vectorization of the forecasts by supporting all forecasts on a common mesh. The common mesh for forecasts may also serve as a mesh to support the next analysis by incorporating the observations, their locations, and any other important factors.

\subsection{Moving Mesh Methods}



Consider a time-dependent PDE written as 
$$\frac{\partial u}{\partial t} = F(u)$$
in an appropriate function space.
The moving mesh approach couples this equation with one that evolves the spatial mesh. The approach we consider is a quasi-Lagrangian approach in which a monitor function or metric tensor, instead of the PDE variable $u$ dieectly, is employed to control mesh concentration. This gives added flexibility in a DA context where it is natural to consider factors beyond the accurate, efficient solution of PDEs used for forecasts. This includes approximation of observational data and potentially targeting regions in space and time for higher resolution, minimizing sensitivities and enhancing DA skill.

In a PDE context, the approach is based (see \cite{Huang2010}) upon an equidistribution condition  that requires all elements to have the same volume
and an alignment condition that requires the elements to be similar to the reference element,
both with respect to a prescribed metric tensor $\M$.
Given a monitor function $\M = \M(x)$, 
these conditions are used to 
define an energy functional 
(see \cite{Huang2010, KHMVVZ22} for details). The mesh equation has the flexibility of employing a potentially arbitrary metric tensor in place of 
a metric tensor derived solely from $u$,
\begin{equation}
    x_t = 
   {\mathcal G}(x,\M). \label{eq:mesheqQL}
\end{equation}
It has been shown that if a mesh begins as nonsingular (that is, the elements have positive volume), it will remain nonsingular for all time under the MMPDE method. Moreover,  the altitudes and volumes of the triangular mesh elements stay bounded below by positive constants that depend only on the metric tensor $\M$, the number of elements, and the initial mesh \cite[Theorem 4.1]{Huang2018}.
In this way, several factors, in addition to those for the PDE solution, may be incorporated into the metric tensor employed to increase the accuracy of observation operators, minimize sensitivities, and enhance DA performance.

\subsection{Combining Meshes}
In order to obtain meshes designed explicitly for DA applications, we combine meshes using their associated metric tensors or mesh density functions. This is done using the so-called metric tensor intersection in (\ref{eq:intersection}). The mesh density function $\M$ obtained in this way is then employed to determine the associated mesh $\X$ using (\ref{eq:mesheqQL}). We form meshes common to all ensemble members designed to support the analysis \cite{KHMVVZ22} or in this work to support the ensemble forecasts and the subsequent analysis. The common mesh is a combination of meshes associated with individual ensemble members, observation locations and their relative impact, and any other factors potentially important to the success of the DA. We adopt the 
philosophy of forming a common mesh once per forecast/analysis cycle.

When remeshing before the analysis, the technique developed in \cite{KHMVVZ22} assumes a forecast ensemble $u^{(i)},\, i=1,...,N_e,$ supported on meshes $(X^{(i)},\M^{(i)})$. Using metric tensor intersection, we form a common mesh by combining the ensemble meshes (through their mesh density functions) with a mesh density function associated with the observations
\[
\M_{\tt ens} = \bigcap_{i=1}^{N_e} \M^{(i)},\,\,\, \M = \M_{\tt ens} \cap \M_{\tt obs}.
\]  
We then form the associated common mesh $\X$, 
interpolate the forecast ensemble to the common mesh, and update analysis. The subsequent ensemble forecasts are formed with independent, time evolving  ensemble meshes $(X^{(i)},\M^{(i)}),\,\, i=1,...,N_e$.

We combine mesh density functions instead of combining meshes directly. Mesh density functions are defined pointwise as $d \times d$ symmetric and positive definite matrices (SPD) and are used to controlled the size, shape, and orientation of elements. To combine mesh density functions we employ the metric tensor intersection.
In particular,
if $A$ and $B$ are SPD matrices of the same size,
then their intersection is denoted by ``$\cap$'' and defined as
\begin{equation}\label{eq:intersection}
A\cap B = P^{-1} \text{diag}\left(\max(1, b_1), \cdots,\max(1, b_d)\right) P^{-T},
\end{equation}
where $P$ is a nonsingular matrix simultaneously diagonalizing $A$ and $B$, i.e., $P A P^T = I$ and $P B P^T = \text{diag}(b_1,\cdots, b_d)$.
For two mesh density functions $A$ and $B$ (meaning they are SPD-matrix-valued functions), the mesh associated with $A\cap B$ will combine the higher resolution of the meshes associated with $A$ and $B$.

The strategy developed in this work is based on forming mesh common to all ensemble members after the analysis is performed to prepare for the next set of ensemble forecasts and the subsequent analysis.
Given analysis ensemble $u^{(i)},\, i=1,...,N_e,$ supported on common mesh $(\X,\M)$, we employ
a (random) subset of $(N_{se})$ of the analysis ensemble to perform pre-forecasts on independent, time-evolving meshes. During the pre-forecasts, accumulated (over time) mesh density functions $\M_{\tt accum}^{(i)},\, i=1,...,N_{se}$ are formed.
A robust look ahead mesh $(\X,\M)$ is formed based on accumulated mesh density functions and the subsequent observations for
next analysis: 
\[
\M_{\tt ens} = \bigcap_{i=1}^{N_{se}}\M_{\tt accum}^{(i)},\,\,\, \M = \M_{\tt ens} \cap \M_{\tt obs}.
\]
The look ahead mesh $(\X,\M)$ is employed as a common mesh to form ensemble forecasts and support the subsequent analysis. Other than intermediate meshes used to form the look ahead mesh, a single nonuniform mesh is employed throughout each forecast/analysis cycle which decreases error due to interpolation between meshes and potentially improves the computationally complexity of forming the ensemble forecasts.

Mesh density functions $\M = \M(\bx)$ are employed to control the size, shape, and orientation of the mesh being generated. Metric tensors have been developed based on a number of different criteria, see \cite{Huang2003} for a survey of choices.

\section{Robust meshes}\label{RobustMeshes}

Developing meshes that are robust over a range of similar behaviors provides advantages in terms of efficiency (improved vectorization, less need for remeshing) and accuracy (by minimizing interpolation). This was done in \cite{KHMVVZ22} by using metric tensor intersection to define meshes that are an ``average'' of the meshes for all the ensemble members
\begin{equation}
\M^{ens} = \M^{1}\cap\M^{2}\cap...\cap\M^{N_e}. 
\end{equation}
A similar idea can be applied to obtain meshes/metric tensors that remain fixed over longer time intervals by intersecting relevant temporal meshes, e.g., over an  observation time scale $\Delta t$, to form 
\begin{equation}
\M^{[0,\Delta t]} = \M^{t_0}\cap\M^{t_1}\cap...\cap\M^{t_M},
\end{equation}
where $t_j = \alpha_j\Delta t$, $j=1,...,M$, for 
$0\leq \alpha_0 < \alpha_1 <\cdots < \alpha_M\leq 1$.

The metric tensors can be formed when predictions are made using the mean and/or (small) ensemble forecasts by evolving from the metric tensor corresponding to the current mesh.
Metric tensor intersection is not the only means available to enhance existing meshes. This may also be done by signaling the need for higher resolution in different regions of the spatial domain and then performing mesh refinement and/or rezoning in these regions.

\subsection{Robust Look Ahead Meshes}\label{LAHsec}

In this section we develop robust look ahead meshes that will be used as both a common mesh for all forecasts and for the subsequent analysis. This contrasts with the approach in \cite{KHMVVZ22} in which a common mesh was formed before the analysis and then the individual ensemble meshes were allowed to evolve independently during the subsequent ensemble forecasts. To develop look ahead meshes, we employ the analysis mean or a subset of the analysis ensemble solutions to determine a mesh that is robust in time and across ensemble members. We fix that mesh over the entire forecast time and for all ensemble members. We illustrate how this mesh is formed using metric tensor intersections.



\subsection{Common Meshes for Forecasts and Analysis}

We next develop robust look ahead meshes. We seek to form a mesh that is effective for all ensemble members, fixed over the time between observations for the ensemble forecasts, and designed to minimize error in the subsequent analysis update. To do so we form a ``look ahead'' mesh using the previous analysis mean or a subset of the analysis ensemble members as initial condition(s).

We employ operations of the following form.
\begin{itemize}
    \item $\M_{\tt new} = \M({\mathbf U},\X_{\tt old})$; [Solution ${\mathbf U}$ supported on $\X_{\tt old}$],
    \item $\X_{\tt new} = \X([t_j, t_{j+1}],\X_{\tt old},\M)$; [Evolve mesh over $[t_j,t_{j+1}]$ starting from mesh $\X_{\tt old}$].
\end{itemize}
One typically employs $\M_{\tt new}$ as the metric tensor $\M$ in the evolution to the new mesh $\X_{\tt new}.$ We aim to form look ahead meshes for an interval $[t,t+\Delta t]$ between subsequent observation times. Given a discretization of the interval $[t,t+\Delta t]$ with $t_0 = t<t_1< \cdots <  t_J = t+\Delta t,$ we recursively form accumulated metric tensors $\M[t_0,t_{j+1}] := \M[t_0,t_{j}]\cap \M(t_{j+1})$ for $j=0,1,...,J-1.$ The mesh is updated as
\begin{equation}\label{MeshUpdate}
\X(t_{j+1})=\X([t_j,t_{j+1}], \X(t_{j}), \M[t_0,t_{j+1}])
\end{equation}
using the accumulated metric tensor as opposed to the standard choice of the current metric tensor $\M(t_{j+1}).$

To form the look ahead mesh we combine the accumulated metric tensors over $[t,t+\Delta t]$ for each of the members of the small ensemble together with a metric tensor for the observations as the next analysis time $t+\Delta t.$
For $i=1,...,N_{e,s}$, let $\M^{(i)}[t,t+\Delta t]$ denote the accumulated metric tensor for the $i$th ensemble member
and then we form either 

 \begin{equation}\label{LAHMesh2}
\M_{\tt LAH}[t,t+\Delta t]:=\bigcap_{i=1}^{N_{se}} \left[\M^{(i)}[t,t+\Delta t] \cap \M^{\tt obs}(t+\Delta t)\right],
\end{equation}
or
\begin{equation}\label{LAHMesh1}
\M_{\tt LAH}[t,t+\Delta t]:=\bigcap_{i=1}^{N_{se}} \left[ \bigcap_{t_j\in [t,t+\Delta t]} \M^{(i)}[t,t_j]\cap \M^{\tt obs}(t_j)\right] ,
\end{equation}
or with $\M^{\tt obs}(t+\Delta t)$ replaced with $\M^{\tt obs}[t,t+\Delta t]$.
After smoothing the metric tensor $\M_{\tt LAH}$, the look ahead mesh is obtained as
\[
\X_{\tt LAH}[t,t+\Delta t]=\X([t,t+ \Delta t],\X(t),\M_{\tt LAH}[t,t+\Delta t]).
\]
This mesh incorporates variations between ensemble members over time with a focus on the next set of observational data.



\subsection{Robust Meshing Example: Traveling Waves in Reaction-Diffusion}

\begin{figure}
    \centering
    \includegraphics[width=3.0true in, trim= 40 270 20 270,clip]{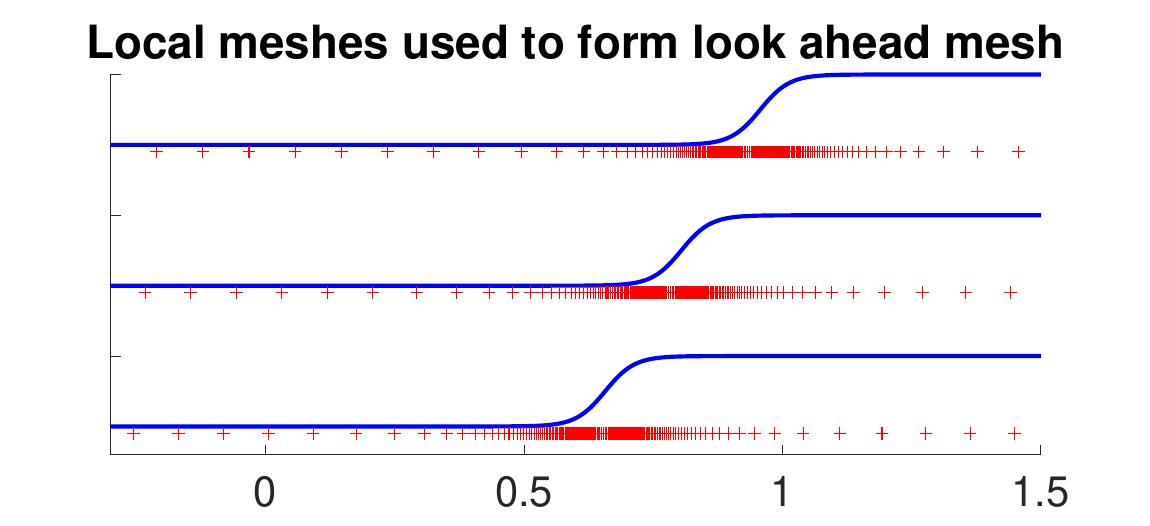}
    \includegraphics[width=3.0true in, trim= 40 270 20 270,clip]{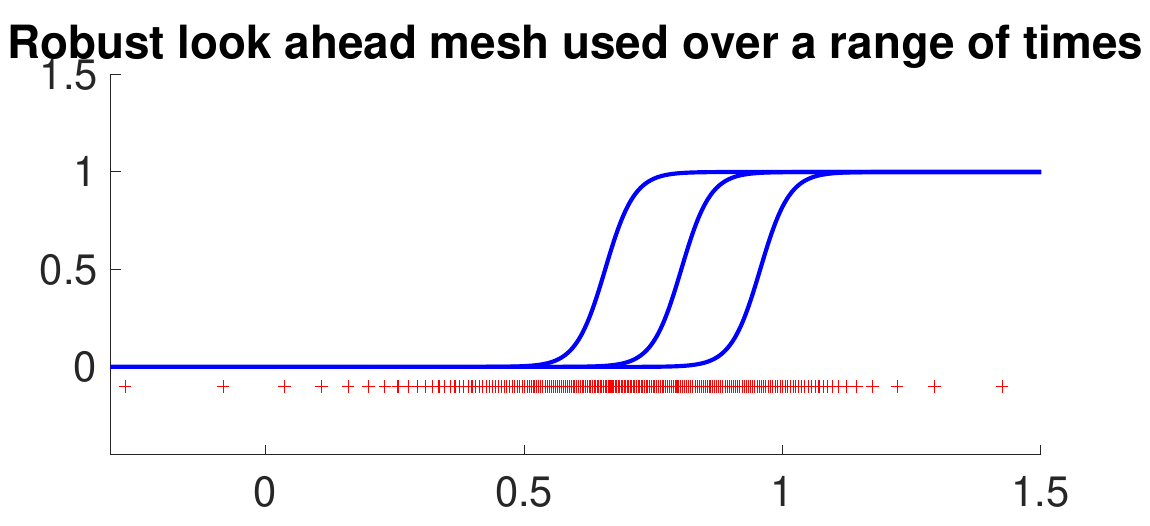}
    \caption{Combining meshes to form a robust, commmon, look ahead mesh.}
    \label{fig:LAH1}
\end{figure}

To illustrate some of the ideas behind the construction of robust meshes, e.g., ensemble and look ahead mesh, 
we consider a bistable reaction-diffusion equation, the Nagumo equation, that is known to possess traveling wave solutions. For this equation with a cubic nonlinearity the existence of traveling wave type solutions has been established for a moving mesh method with an arc length monitor function/mesh density function/metric tensor, see \cite{HupkesEVVMMRDIII}. The specific form of the Nagumo PDE we consider is 

\begin{equation}\label{RDE}
    u_t = \epsilon^2u_{xx}-f(u;a),\,\,\, 0<x<L,\,\,\, t>0
\end{equation}
together with Dirichlet boundary conditions using the exact traveling wave solution $u(0,t)=u_{\tt exact}(0,t),\,\,\, u(L,t)=u_{\tt exact}(L,t).$ We discretize the problem using finite differences on a nonuniform mesh and the mesh evolves using an arc-length monitor function (mesh density function), see \cite{HupkesEVVMMRDIII,KHMVVZ22}. For the cubic nonlinearity $f(u)=u(u-1)(u-a)$, the exact from of the wave profile $\phi$ and wave speed $c$ are known, 
\begin{equation}\label{TWRDE}
    u(x,t) = \phi(x-ct) = \frac{1}{2}\left[1 + \tanh(\frac{x-ct}{2\sqrt{2}\epsilon})\right],\,\,\, a = \frac{1}{2} + \frac{c}{\epsilon \sqrt{2}}.
\end{equation}
In Figure \ref{fig:LAH1} we illustrate combining local meshes that change with time to form a single mesh that is employed over an interval of time. The left figure shows the wave form and associated mesh at three times and then in the right figure the combined mesh that supports the accurate, efficient propagation of the traveling wave over the entire time interval.

%
%
%

\begin{figure}
    \centering
    \subfloat[$N_{se} = 2$ and $\bP_0 = 0.1 \cdot I$]{\includegraphics[width=3.0true in, trim= 35 160 20 180,clip]{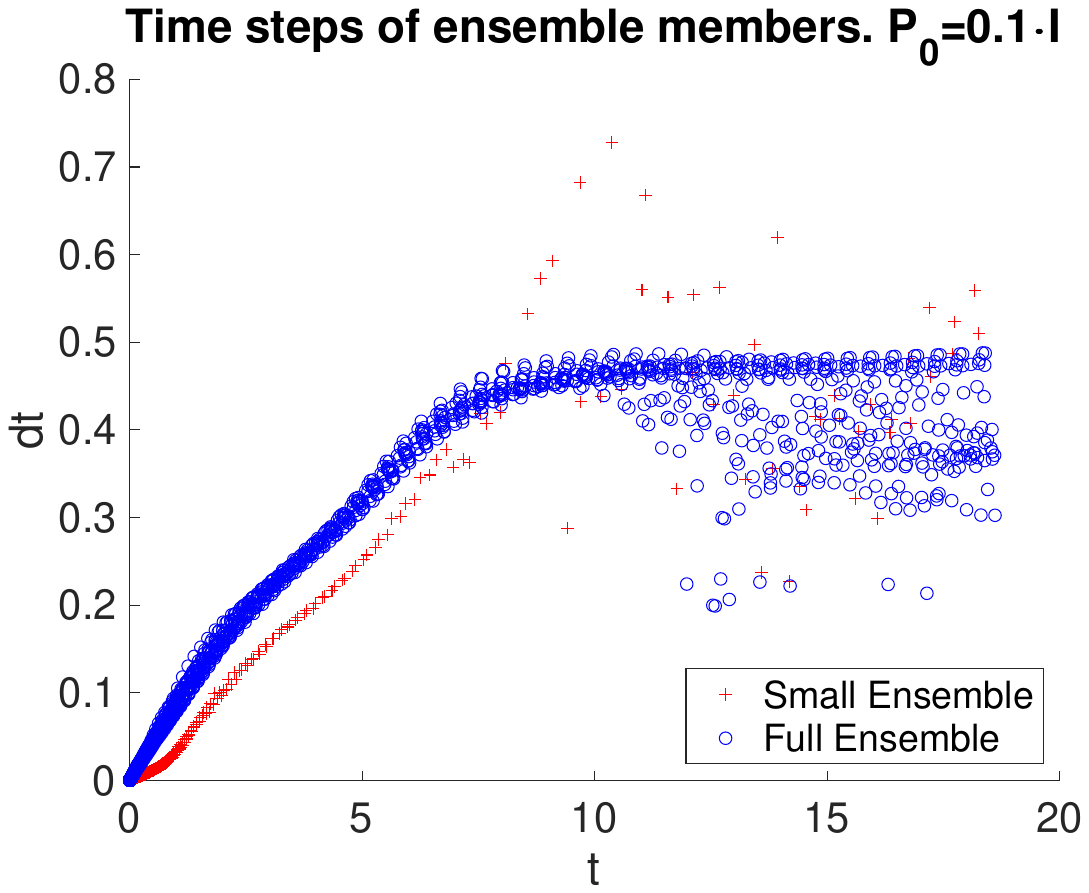}}
    \subfloat[$N_{se} = 2$ and $\bP_0 = 0.01 \cdot I$]{\includegraphics[width=3.0true in, trim= 35 160 20 180,clip]{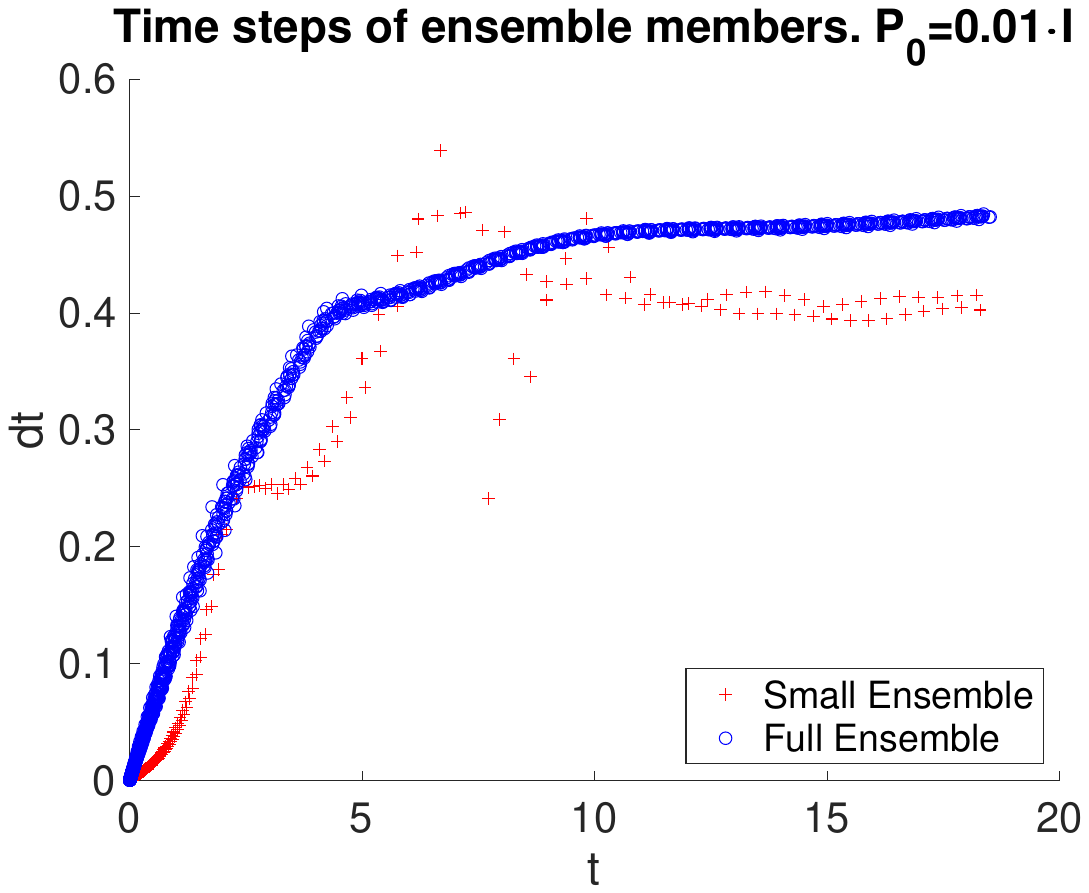}}
    \\
    \subfloat[$N_{se} = 8$ and $\bP_0 = 0.1 \cdot I$]{\includegraphics[width=3.0true in, trim= 35 160 20 180,clip]{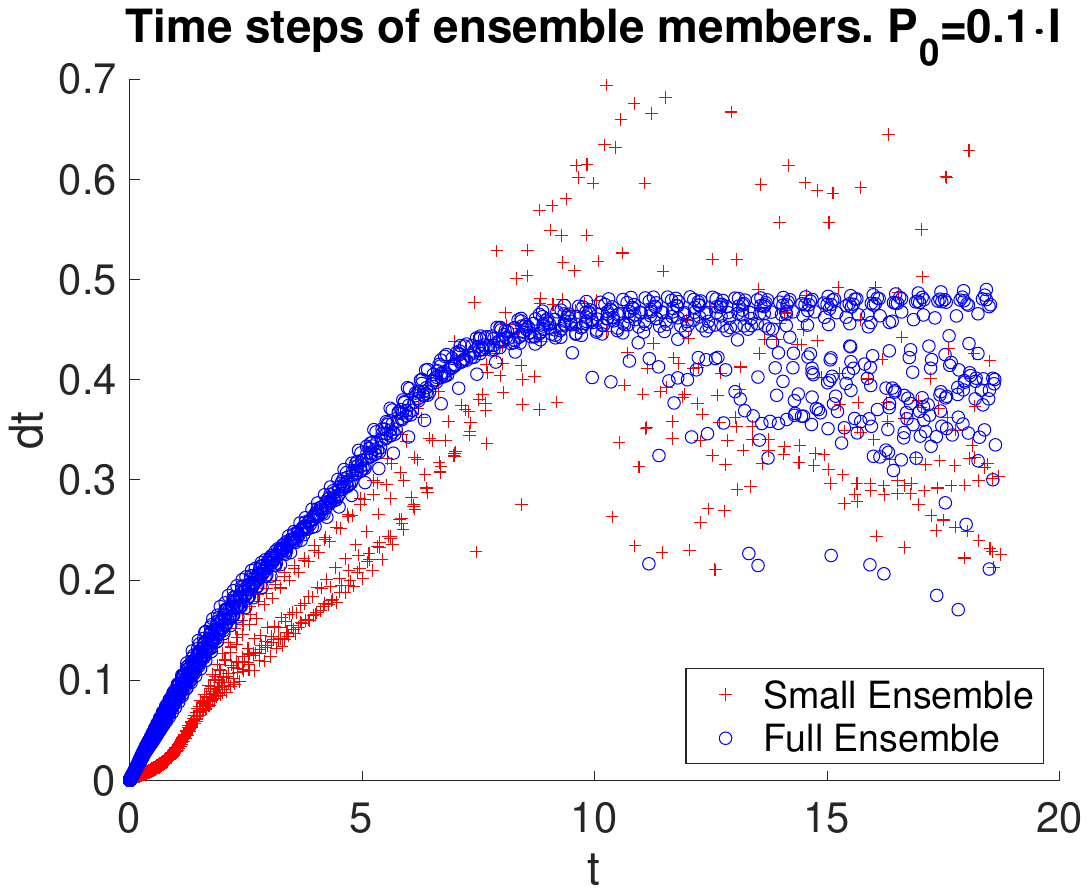}}
    \subfloat[$N_{se} = 8$ and $\bP_0 = 0.01 \cdot I$]{\includegraphics[width=3.0true in, trim= 35 160 20 180,clip]{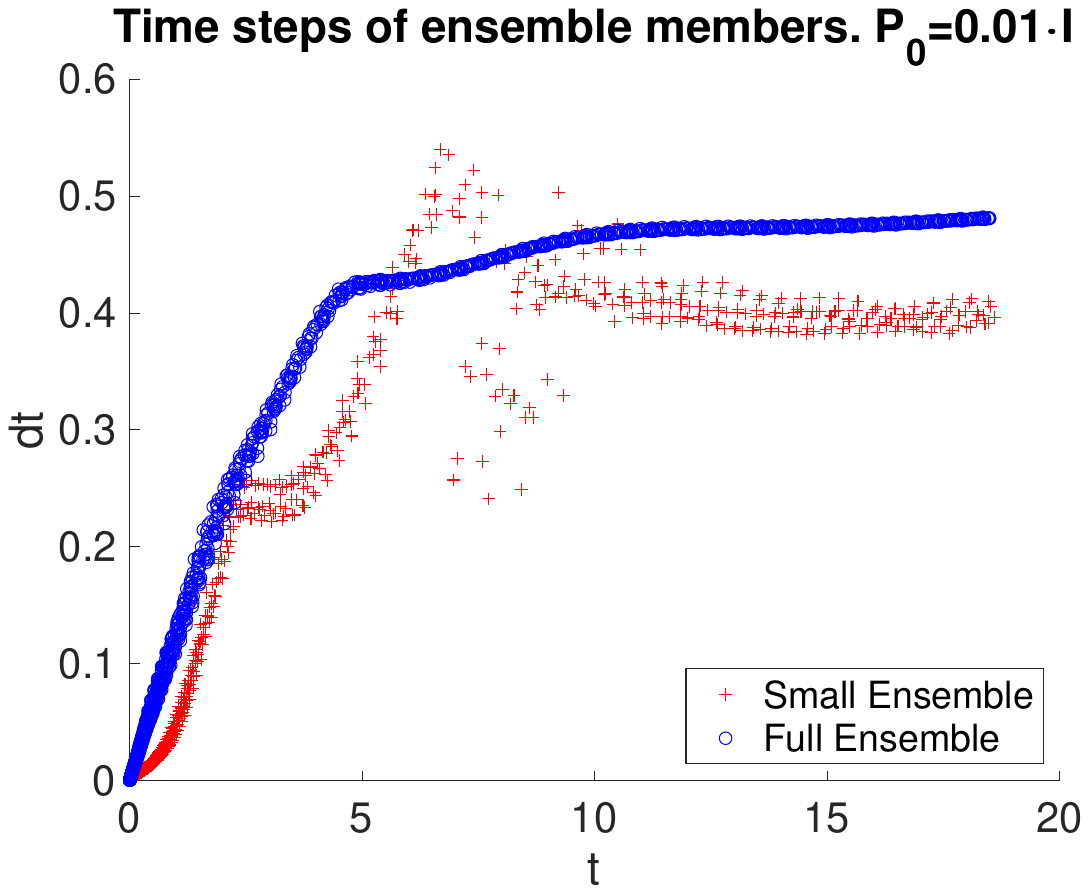}}
    \caption{Timesteps of the small ensemble employed when forming the LAH  mesh and timesteps of the full ensemble when using the LAH as a fixed mesh. $\bP_0$ is the covariance matrix of the initial perturbations of the ensemble members.}
    \label{fig:LAH2}
\end{figure}

For Figure \ref{fig:LAH2} we consider the differential equation (\ref{RDE}) with parameters $\epsilon^2 = 10^{-2}$, $a = 0.95,$ and wave speed $c$ calculated as in (\ref{TWRDE}). We calculate the look ahead mesh using small ensemble sizes $N_{se}=2,8$ in which the small ensemble is created by perturbing the traveling wave initial condition by drawing from a Gaussian distribution with mean zero and covariance $\bP_0 = \alpha\cdot I$ for $\alpha=0.1, 0.01.$ 
Each member of the small ensemble solutions evolves on its own time dependent mesh independent of the other members. The starting mesh for all ensemble members is taken to be the mesh of the unperturbed initial condition. Members of the full ensemble are then evolved using the LAH mesh formed by the accumulated mesh density function over time and over the small ensemble using (\ref{LAHMesh2}) with $M_{\tt obs}=I$. 

We use the time integration scheme
{\tt ode15s} of Matlab with adaptive time stepping
and absolute and relative error tolerances ${\tt tol}_{abs}=10^{-10}$, ${\tt tol}_{rel}=10^{-12}$.
These small tolerances were employed to avoid reaching the maximum stepsize.
In the plots on the left in Figure \ref{fig:LAH2}, with larger variance, the time steps of some ensemble members decrease after $t\approx 10.$ This occurs for both the small ensemble and the full ensemble. In comparison, the plots on the right side in Figure \ref{fig:LAH2} with smaller variance generally maintain much less variation in the time steps between ensemble members. When $P_0 = 0.01\cdot I$ with either $N_{se}=2$ or $N_{se}=8$ a single time step sequence (uniform to all ensemble members, and potentially highly parallel) would likely be effective.

%


\section{Ensemble Based Data Assimilation Techniques}\label{EnsembleDA}


The techniques developed here are applicable to ensemble based DA techniques, including ensemble Kalman filter (EnKF) techniques, particle filters (PFs), as well as hybrid variational techniques, see, e.g., \cite{LSZ, ReichCotter15, BocquetBook}.
In this work we will focus on an EnKF technique called the Local Ensemble Transform Kalman Filter (LETKF) \cite{Hunt} and the use of domain localization using a exponential localization function derived from the mesh density functions. We consider nonlocal observations based upon integrating against the exponential kernel function and derive a mesh density function for these types of observations.

Ensemble DA procedures use an ensemble of solutions to make predictions for the physical state via a two-step process. First, the prediction step uses the physical model in \eqref{DAsetup} to integrate the ensemble members $\{\bu_n^{i}\}_{i=1}^{N_e}$ to make the ensemble forecasts $\{\hat \bu_{n+1}^{i}\}_{i=1}^{N_e}$ where $N_e$ is the number of ensemble members. Second, the analysis step incorporates the observations at time $t_{n+1}$ to adjust the prediction $\{\bu_{n+1}^{i}\}_{i=1}^{N_e}$. The Ensemble Kalman Filter (EnKF) does this through the following sequential process:
\begin{equation}
    \text{Forecast:} \quad \begin{cases}\hat \bu_{n+1}^{i} &= \Psi(\bu_n^{i})+\xi_n^{i},\quad i = 1, \dots, N_e\\
    \hat \m_{n+1} &= \frac{1}{N_e}\sum_{i=1}^{N_e} \hat \bu_{n+1}^{i}, \\
    \bP^f_{n+1} &= \frac{1}{N_e-1}\sum_{i=1}^{N_e}\left(\hat \bu_{n+1}^{i}-\hat \m_{n+1}\right)\left(\hat \bu_{n+1}^{i}-\hat \m_{n+1}\right)^T ,
    \end{cases} \label{eq:DA_pred}
\end{equation}
\begin{equation}
    \text{Analysis:} \quad \begin{cases} \bK_{n+1} &= \bP^f_{n+1}\bH^T\left(\bH\bP^f_{n+1} \bH^T + \bR\right)^{-1}, \\
    \bu_{n+1}^{i} &= \hat \bu_{n+1}^{i} + \bK_{n+1}[\by_{n+1} - \mathcal{H}(\hat \bu_{n+1}^{i})], \quad i = 1, \dots, N_e
    \end{cases}\label{eq:DA_analysis}
\end{equation}
where  $\mathbf{H}$ is the linearization of $\mathcal{H}$, and $\mathbf{K}_{n+1}$ is the Kalman gain matrix.

In this work we specifically employ a square root filter with localization LETKF which determines the analysis update \eqref{eq:DA_analysis} as a linear combination of the ensemble solutions. With an ETKF technique, the sample forecast covariance matrix $\bP_{n+1}^f$ in \eqref{eq:DA_pred} is written as $\bP_{n+1}^f = \bX_{n+1}\bX_{n+1}^T$ and the analysis update in \eqref{eq:DA_analysis} is written for $i = 1, \dots, N_e$ as 
\begin{equation}
    \bu_{n+1}^{i} = \hat \bu_{n+1}^{i} + \bX_{n+1}\mathbf{w}_{n+1}^i,\,\,
    \mathbf{w}_{n+1}^i= \bX_{n+1}^T\bH^T\left(\bH\bP^f_{n+1} \bH^T + \bR\right)^{-1}\bI_{n+1}^i, 
    \label{eq:DA_analysisETKF}
\end{equation}
where $\bI_{n+1}^i = \by_{n+1} - \mathcal{H}(\hat \bu_{n+1}^{i})$ is the innovation for the $i$th ensemble member. This may be further modified by application of the Sherman-Morrison-Woodbury formula and is commonly written as (see \cite{BocquetBook})
\begin{equation}\label{ETKFupdate}
\mathbf{w}_{n+1}^i = \bY_{n+1}^T (\bY_{n+1} \bY_{n+1}^T + \bR)^{-1}\bI_{n+1}^i
                   = (\bI + \bY_{n+1}^T\bR^{-1}\bY_{n+1})^{-1}\bY_{n+1}^T\bR^{-1}\bI_{n+1}^i
\end{equation}
where $\bY_{n+1} = \bH\bX_{n+1}.$

We will employ a simple form of domain localization employed in \cite{Hunt} based upon a user-defined localization radius for each spatial point in the model domain. This results in local analysis updates that can be formulated in terms of simple projections of the form $\bP^{(k)} = \bS^{(k)} (\bS^{(k)})^T$ and $\bP^{(k)}_o = \bS^{(k)}_o (\bS^{(k)}_o)^T$ where for each spatial point $\bx_k$ with localization radius $\br_k$, $\bP^{(k)}$ is the orthogonal projection onto the spatial points $\bx_j$ such that $\|\bx_k - \bx_j\|\leq \br_k$ and $\bP^{(k)}_o$ is the orthogonal projection onto the the observation locations $\bx_j^o$ such that $\|\bx_k - \bx_j^o\|\leq \br_k$. In this way the form of a local update corresponding to the spatial location $\bx_k$ is obtained from (\ref{eq:DA_analysisETKF}), (\ref{ETKFupdate}) by replacing 
(ignoring subscripts and superscripts)
$\bX\to \bS^{(k)}(\bS^{(k)})^T\bX$, $\bY \to \bH\bS^{(k)}(\bS^{(k)})^T\bX$, and $\bR\to \bS^{(k)}_o(\bS^{(k)}_o)^T \bR \bS^{(k)}_o(\bS^{(k)}_o)^T$. This results in local updates of the form 
\begin{equation}\label{LETKFupdate}
\bS^{(k)}(\bS^{(k)})^T\bu_{n+1}^i = \bS^{(k)}(\bS^{(k)})^T\hat \bu_{n+1}^{i} + \bS^{(k)}(\bS^{(k)})^T\bX_{n+1}(\mathbf{w}_{\bS^{(k)}})_{n+1}^i,
\end{equation}
where 
\begin{align}
\label{LETKFupdatewS}
(\mathbf{w}_{\bS^{(k)}})_{n+1}^i & = (\bI + \bY_{S^{(k)}}^T\bS^{(k)}_o((\bS^{(k)}_o)^T\bR\bS^{(k)}_o)^{-1}(\bS^{(k)}_o)^T\bY_{S^{(k)}})^{-1}
\\
& \qquad \cdot \bY_{S^{(k)}}^T\bS^{(k)}_o((\bS^{(k)}_o)^T\bR\bS^{(k)}_o)^{-1}\cdot [(\bS^{(k)}_o)^T \bI_{n+1}^i],
\notag
\end{align}
where $\bY_{S^{(k)}} = \bH \bS^{(k)} (\bS^{(k)})^T\bX_{n+1}.$ The analysis update is obtained from the local analysis updates using (\ref{LETKFupdate}) as
\begin{equation}\label{LETKFanalysisupdate}
 \bu_{n+1}^i = \sum_k \be_k \be_k^T[\bS^{(k)}(\bS^{(k)})^T\bu_{n+1}^i], 
\end{equation}
where $\be_k \be_k^T$  is the projection onto the $k$th solution component, $\be_k^T \bu_{n+1}^i = \bu_{n+1}^i(\bx_k)$ and $\sum_k \be_k\be_k^T = \bI$.

For multi-component systems, there are some options on how to perform the local analysis update. Simplest among these is to perform the local analysis update separately for each component using (\ref{LETKFupdate}), (\ref{LETKFupdatewS}), (\ref{LETKFanalysisupdate}). Alternatively, the analysis update may be performed simultaneously on all components. This has the advantage that the covariance information between model components is included, and the update is straightforward to facilitate in our framework since all model components are supported on the same mesh. The coupled approach is obtained from  (\ref{LETKFupdate})--(\ref{LETKFanalysisupdate}) by considering $\bu,\X,$ etc. as multi-component variables and then replacing $\bS\bS^T$ and $\bS_o\bS_o^T$ by the block diagonal matrices ${\tt diag}(\bS\bS^T,...,\bS\bS^T)$ and ${\tt diag}(\bS_o\bS_o^T,...,\bS_o\bS_o^T).$



\subsection{Localization and Local Analysis for nonlocal Observation Operators}

Localization is performed via a local analysis, domain localization, in which the analysis update at a given spatial location is restricted to observations that lie within a prescribed neighborhood of this spatial location. This allows for an analysis update that may be performed in parallel as independent updates depending on spatial location. The use of localization is important to the success of DA techniques, and there is a well established theory \cite{MorzfeldHodyss} for local, pointwise defined observations. For nonlocal observations, the use of localization techniques is still imperative but the techniques are more ad hoc \cite{Fertignonlocal, PeterJanFAMSnonlocal}.

Following \cite{KHMVVZ22} we employ an exponential localization based on the determinant of the local metric tensors (MT localization). This is used for  domain localization scheme with LETKF. At each observation time the localization radius is calculated for each node:
\begin{equation}\label{eq:dynloc}
r_i = r_0\;  e^{-\frac{(d_i-d_{min})}{2 d_{min}}} , \quad i = 1, ..., N_u,
\end{equation}
where $d_i = \min(\det(\M_K^m(\bx_i)),c)$, $d_{min} = \min_i d_i$, and $c \geq d_{min}$ and $r_0 > 0$ are the parameters. In this way the localization radius for all points is bounded above and below as 
\[
r_0e^{-\frac{(c-d_{min})}{2 d_{min}}}\leq r_i\leq r_0. 
\]
We determine the cutoff value $c$ in our experiments as $c=(d_{max}+d_{min})/2$ which is recomputed before the current analysis update.
This may be employed to define a neighborhood $N(\bx_i;I)$ for the local analysis as a Euclidean ball of radius $r_i$ centered at the mesh point $\bx_i$ so the local analysis is performed using observations associated with locations $\hat \bx_j$ with $\hat \bx_j \in N(\bx_i;r_i,I).$



When employing such nonlocal observation operators (as opposed to pointwise defined observation operators) with LETKF using domain localization, there is not a well established way \cite{Fertignonlocal, PeterJanFAMSnonlocal} of defining the localization.
Given a point $\bx\in\Omega$ with a localization function defining a neighborhood $N(\bx)$ about $\bx$, a nonlocal observation centered at $\hat \bx_i$ will be used to perform the update at the point $\bx$  provided $\hat \bx_i \in N(\bx)$. In addition, we include the entire support of the observation centered at $\hat \bx_i$ which effectively increases the size of the neighborhood $N(\bx)$ about $\hat \bx_i.$

\begin{remark}
An alternate strategy that we have not pursued here is to increase the domain of localization for nonlocal observation operator by only increasing the localization to include those mesh points in the support of the nonlocal observations as opposed to increasing the entire radius.
    
Also note that the norm induced by the local mesh density function/metric tensor may be used to define a potentially anisotropic neighborhood for the local analysis.  This changes the neighborhood to $N(\bx_i;r_i,\M_i)=\{\bx|(\bx-\bx_i)^T\M_i(\bx-\bx_i)^T < r_i\}.$
\end{remark}


\subsection{Mesh Density Functions}\label{MTsec}


There is a wealth of knowledge on the use of metric tensors/mesh density functions for the efficient, accurate numerical solution of both time dependent and time independent PDEs. We want to adapt some of these well-developed ideas for application in DA contexts. Foremost among these is the development of mesh density functions for observational operators and data. We note that mesh density functions for DA can be defined for scalar-valued functions and by extension vector-valued functions with sufficient smoothness. 


In this work, for the model forecast solutions, we employ a Hessian-based metric tensor defined for each element $K \in \T_h$ as
\begin{equation}\label{mer}
   \M_{K} =\det \left(I+\frac{1}{\alpha_h}|H_K(u)|
        \right)^{-\frac{1}{d+4}}
   \left(I+\frac{1}{\alpha_h}|H_K(u)|\right) ,
\end{equation}
where $H_K(u)$ is the Hessian or a recovered Hessian of the state vector $u \in \R^d$ on the element $K$; $|H_K(u)| = Q
\text{diag}(|\lambda_1|,...,|\lambda_d|)Q^T$ with $Q \text{diag}(\lambda_1,
...,\lambda_d)Q^T$ being the eigen-decomposition of $H_K(u)$; and $\alpha_h$
 is a regularization parameter that may be prescribed or defined optimally using an algebraic equation.
The metric tensor given by equation \eqref{mer} is known to be optimal for the $L^2$-norm of linear interpolation error \cite{Huang2003}. 

In \cite{KHMVVZ22} an ad hoc metric tensor was defined for observations supported at points within the spatial domain. This was based upon the use of Gaussian type indicator functions that depend only on the spatial location of the observation and not the corresponding value of the ensemble solutions in a neighborhood of the observation locations.
To avoid fixing the observation location as a point in the mesh, the metric tensor for the observation locations concentrates mesh points
near the observation locations.  The metric tensor for the location of point observations employed in \cite{KHMVVZ22} (see also\cite{DHHLY18}) has the form ($\sigma=1/2$ was employed)
\begin{equation}\label{MobsAdHoc}
 \M_{obs} = I + \sum_{j=1}^{N_y}\chi\left(\|\bx-\bx_j^o(t)\|\right)I,\,\,   \chi(w) = \left[e^{w^2/\sigma^2}-1+\frac{2}{\max\limits_{K} \sqrt{\det(\M_K^{ens})}}\right]^{-1}.
\end{equation}


We next develop a general metric tensor for a class of nonlocal observations that depend not only the location of the support of the observation but also on the corresponding value of the ensemble solutions. 

\subsubsection{Goal-Oriented Metric Tensors for nonlocal Observations}\label{sec:NonLocObs}

Consider $m$ observation locations $\hat \bx_i,$ $i=1,...,m,$ and a kernel function $G=G(\bx),$ $\bx\in\Omega,$ and (linear) nonlocal observation operators $H$
of the form 
\begin{equation}
    Hu(\hat \bx) = \int_\Omega u(\bx)G(\hat \bx - \bx)d \bx
\end{equation}
where $u(\bx)$ is the physical PDE variable at a fixed time $t$.
The observational data we consider is of the form
\begin{equation}
    \tilde y_i = (Hu^t)(\hat \bx_i) = \int_\Omega u^t(\bx) G(\hat \bx_i-\bx)d \bx,\,\, i=1,...,m,
\end{equation}
where $u^t$ represents the unknown truth. We approximate $(Hu)(\hat \bx_i)$ as
\begin{equation}\label{HuApprox}
     Hu(\hat \bx_i) = \int_\Omega u(\bx)G(\hat \bx_i - \bx)d \bx
     \approx \sum_{K} |K|G(\hat \bx_i - \bx_K)\cdot \frac{1}{d+1}\sum_{j=0}^d u(\bx_j^K)
\end{equation}
using the average value of $u$ over the $d+1$ vertices in the current element. In our numerical experiments we will employ the Gaussian kernel function
\begin{equation}
  G(\bx) = \frac{1}{(\sqrt{2\pi}\delta)^d} e^{-\frac{\|\bx\|^2}{2\delta^2}}
\end{equation}
where $\delta > 0$ is a scaling parameter.

Obtain the metric tensor $\M_K$ for the element $K,$
\begin{equation}
    \M_K = \det(A_K)^{-\frac{1}{d+2}} A_K,\quad
    A_K = I + \frac{1}{\alpha_h}|H_K|\sum_{i=1}^{N_y} |G(\hat \bx_i - \bx_K)|
\end{equation}
where $d$ is the spatial dimension, $H_K$ is the approximate Hessian for $u$ on the element $K$, $\bx_K$ is the centroid of $K$, and $\alpha_h$ is a regularization parameter than may be set, e.g., $\alpha_h = 1,$ or an optimal value employed. See Appendix \ref{AppA} for further details of the derivation of the mesh density function for these nonlocal observations.

In our implementation we replace the globally defined (over the entire domain $\Omega$) convolution based observation operator \eqref{HuApprox}
with
\begin{equation}\label{HuApproxLoc}
     (Hu)(\hat \bx_i) = \int_{B_{r_i}(\hat \bx_i)} u(\bx)G(\hat \bx_i - \bx)d \bx
     \approx \sum_{\{K:\bx_K\in B_{r_i}(\hat \bx_i)\}} |K|G(\hat \bx_i - \bx_K)\cdot \frac{1}{d+1}\sum_{j=0}^d u(\bx_j^K),
\end{equation}
where $B_{r_i}(\hat \bx_i)$ denotes a ball of radius $r_i$ about $\hat \bx_i$.


\section{Computational Techniques and DA Parameters}\label{DACompParams}

In our numerical experiments we consider three model problems and a number of different observation scenarios. The first two problems are hyperbolic PDEs in two space dimensions. We outline the numerical technique employed for these problems that is based upon a discontinuous Galerkin approximation on nonuniform meshes in space and variable time stepping using strong stability preserving time stepping with embedded pairs. The third example is discretized using nonuniform finite differences in space and fixed step size linearly implicit backward Euler time stepping. In addition to the discretization of these PDEs, we describe a number of parameters that will be used for the data assimilation with LETKF.

\subsection{Moving Mesh DG, Adaptive Time Stepping, Meshing Functionals, and Interpolation}

For our experiments, we employ
the discontinuous Galerkin (DG) method (see \cite{DGforDA} for recent work on tailoring DG to DA techniques) for spatial discretization together with time stepping using embedded pairs of strongly stability preserving (SSP) Runge-Kutta schemes.  The DG method has upwinding naturally built into its discretization (via the choice of numerical fluxes and limiters) and results in stable computation for hyperbolic problems. DG is a type of finite element method with the trial and test function spaces consisting of discontinuous, piecewise polynomials. It is known to be a particularly powerful numerical tool for the simulation of hyperbolic problems and has the advantages of high-order accuracy, local conservation, geometric flexibility, suitability for handling mesh-adaptivity, extremely local data structure, high parallel efficiency, and a good theoretical foundation for stability and error estimates. Two types of interpolants are employed primarily to interpolate ensemble members from the previous common mesh to the current common mesh after each analysis update, a linear interpolant and a higher order DG-interpolation scheme \cite{Huang2021} that is conservative and positivity preserving.  

Specifically, we employ the quasi-Lagrangian moving mesh DG method (MMDG), see \cite{Huang2021}.  The resulting method is (mass) conservative so that the model error maintains mean 0. In our computation, we use piecewise linear polynomials (P$^1$-DG) for spatial discretization and we employ the third-order explicit SSP embedded Runge-Kutta pair (SSP-ERK(4,3)) of \cite{FEKETE2022114325} using a standard step-size selection strategy (e.g., see Hairer and Wanner \cite{HW91}). The meshing functional employed is based upon the equidistribution and alignment conditions.

Given a mesh density function or metric tensor $\M$, a mesh is uniform in this metric when all elements have a constant volume and are similar to the reference element in the metric $\M$. These conditions are called the equidistribution and alignment criteria \cite[Section 4.1.1]{Huang2010}. To be specific, let $\M_K$ to be the average of $\M(\bx)$ over the element $K$, that is,
\begin{equation}
    \M_K = \frac{1}{|K|}\int_K\M(\bx)d\bx.
\end{equation}
Then the equidistribution condition is given by
\begin{equation} \label{eq:equidistribution}
    \sqrt{\det (\M_K)}\; |K| = \frac{\sigma_h}{N}, \quad \forall K \in \T_h
\end{equation}
where $\sigma_h = \sum_{K\in \T_h}|K|\sqrt{\det(\M_K)}$ and $N$ is the number of simplexes/elements in $\T_h$. The alignment condition is equivalent to
\begin{equation}\label{eq:alignment}
\frac{1}{d}\tr\left(\left(F_K'\right)^{-1}\M_k^{-1}\left(F_K'\right)^{-T}\right)
=\det\left(\left(F_K'\right)^{-1}\M_K^{-1}\left(F_K'\right)^{-T}\right)^{\frac{1}{d}}, \quad \forall K\in\T_h
\end{equation}
where $\det(\cdot)$ and $\tr (\cdot)$ denote the determinant and trace of a matrix, respectively.
 The equidistribution condition (\ref{eq:equidistribution}) requires all elements to have the same volume
while the alignment condition (\ref{eq:alignment}) requires the elements to be similar to the reference element,
both in the metric $\M$.
Given a user-prescribed monitor function $\M = \M(\bx)$, equations \eqref{eq:equidistribution} and \eqref{eq:alignment} can be used to define an energy functional
\begin{equation}\label{energy}
\begin{split}
I_h
=&\frac{1}{3}\sum_{K\in\mathcal{T}_h}|K|\hbox{det}({\mathbb{M}}_K)^{\frac{1}{2}}\left(
\hbox{tr}((F'_K)^{-1}{\mathbb{M}}^{-1}_K(F'_K)^{-T})\right)^{\frac{3 d}{4}}
\\&
+\frac{1}{3} d^{\frac{3 d}{4}}\sum_{K\in\mathcal{T}_h}|K|\hbox{det}
({\mathbb{M}}_K)^{\frac{1}{2}}
\left (\hbox{det}(F'_K) \hbox{det}({\mathbb{M}}_K)^{\frac{1}{2}} \right )^{-\frac{3}{2}}.
\end{split}
\end{equation}
The functional (\ref{energy}) is optimized using a gradient descent technique based upon solving a gradient descent differential equation (see, e.g., \cite{KHMVVZ22}) to implement the mesh update (\ref{MeshUpdate}).










\subsection{Data Assimilation Parameters and Experiments}

A number of parameters are employed to implement LETKF with nonuniform meshes. These include the prior covariance $P_0^b,$ the covariance for the model error $\bQ,$ and the covariance $\bR$ for the observation error. The number of ensemble members is set to $N_e=20$ in most of our experiments. We 
focus on varying the overall localization radius, the inflation factor, and for nonlocal observations, the values of $\delta,$ the kernel parameter, and $r_{obs}$, the radius that defines the support of the integral observations. We also vary the size of the small ensemble $N_{se}$ used to form the robust look ahead mesh that supports all the ensemble forecasts and the subsequent analysis.

%
%
%

\section{Numerical Experiments}\label{Numerics}


We next outline the numerical experiments for three model problems and several different observation scenarios. For the first example, a two space dimension inviscid Burgers' equation, we consider moving, nonlocal observations of two types. We investigate how the DA results depend on the size of the small ensemble used to form the look ahead mesh and then compare how the DA results depend on the parameters describing the nonlocal observations. For the second example, we consider a two space dimension shallow water model configured for gravity waves. We consider uniform, fixed in time arrays of local (point) and nonlocal observations and illustrate how the DA results depend on the density of the observations. 
The first two examples make use of the full moving mesh DG approach to mesh movement while employing SSP Runge-Kutta time stepping techniques with variable step size using embedded pairs. The third example we consider is a coupled system of two 1D Kuramoto-Sivashinsky (KS) equations. Both local and nonlocal observations are considered in fixed arrays either moving or stationary in time. We investigate the dependence on size of the small ensemble used for the look ahead mesh, the use of coupled versus uncoupled analysis updates, as well as two different techniques for accumulating the mesh density function used to determine the look ahead mesh. The coupled KS equations are solved using nonuniform finite differences and backward Euler using a fixed time step. Mesh are updated in time via interpolation after determining the robust updated mesh.

\subsection{Example 1: 2D Inviscid Burgers' Equation}

Our first example is a 2D inviscid 
Burgers' equation (see also \cite{KHMVVZ22})
\begin{equation}
    \label{eq:Burgers2d}
    u_t + u (u_x + u_y) =0,
\end{equation}
with $\Omega = (-0.5,1)\times (-0.5,1)$ and periodic boundary conditions. We employ the initial condition
\[
u = \exp(-\gamma(x^2 + y^2)),
\]
with $\gamma = -\log(10^{-16})$. The solution will
have a Gaussian bump that propagates diagonally to the upper right corner of the domain and through the use of periodic boundary conditions continues to propagate from the lower left corner to the upper right corner of the domain.

\begin{figure}
    \centering
 \includegraphics[width=3.0true in, trim= 160 270 160 270,clip]{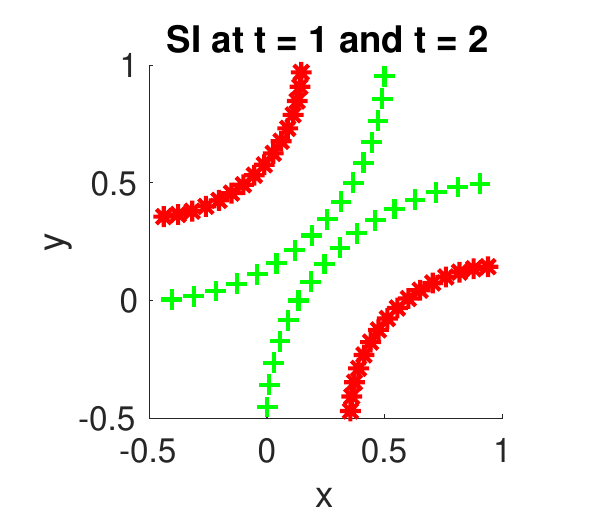}
 \includegraphics[width=3.0true in, trim= 160 270 160 270,clip]{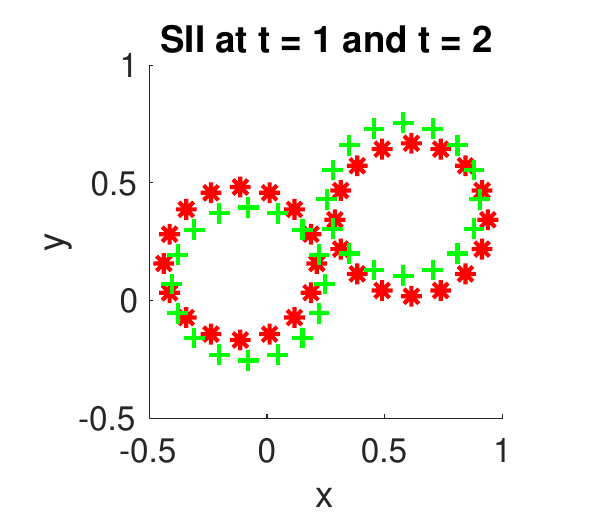}
    \caption{Snapshots of moving observations for observation scenarios I and II. Green denotes current observation locations and red denotes the locations at the previous observation time.}
    \label{fig:BurgersSISII}
\end{figure}

We consider two observations scenarios, both of which are based upon moving arrays of observers. In the first scenario (\ref{TwoQuarter1}), (\ref{TwoQuarter2}) the observations are supported on two quarter circles that moving with a fixed velocity through the domain. In the second scenario (\ref{TwoCircs}) the observers are supported on two circles that rotate about the center point of the domain with a fixed velocity. Plots of the observer locations for these two scenarios with $v=0.25$ are shown in Figure \ref{fig:BurgersSISII}.

{\bf Observation Scenario I: Supported on Two Quarter Circles}

\begin{equation}\label{TwoQuarter1}
\begin{split}
&    (x_1,y_1) + r(t)\cdot (\cos(\theta_j + \pi/2), \sin(\theta_j + \pi/2),
    \\
&    (x_2,y_2) + r(t)\cdot (\cos(\theta_j - \pi/2), \sin(\theta_j - \pi/2),
\end{split}
\end{equation}
with $(x_1,y_1) = (1, -1/2), (x_2,y_2) = (-1/2,1)$ where for some velocity $v$,
\begin{equation}\label{TwoQuarter2}
    r(t) = \frac{1}{2}\sin^2(vt\cdot \pi/2),\quad
    \theta_j \in (0, \pi/2),\,\, j=1,...,N_y/2.
\end{equation}

{\bf Observation Scenario II: Supported on Two Circles Circling a Center Point}
\begin{align}\label{TwoCircs}
\begin{split}
&   (x_0,y_0) + r_0 (\cos(vt), \sin(vt))+ r_1(\cos(\theta_j),\sin(\theta_j)),\\
&   (x_0,y_0) + r_0 (\cos(vt+\pi), \sin(vt+\pi))+ r_1(\cos(\theta_j),\sin(\theta_j)),
\end{split}
\end{align}
for $\theta_j\in [0,2\pi], j=1,...,N_y/2,$ where we take $x_0=y_0 = 0.25$, $r_0 = 0.375$, $r_1 = 0.325.$ 

In our first set of experiments we compare the computational cost of forming the look ahead mesh and the subsequent ensemble forecasts employing the look ahead mesh for all ensemble members with forecasts done in which each ensemble members dynamically determines its own mesh over the forecast cycle.
For the comparison, 
we fix the ensemble size $N_e = 20,$ the initial prior covariance $\bP_0^b = 10^{-2}\cdot I$, the model error covariance and observation error covariances $\bQ = 10^{-1}\cdot I$ and $\bR = 10^{-2}\cdot I$, respectively. We set the covariance inflation factor $\rho = 1.1$ with $N_y = 16$ stationary quarter-circle observations corresponding to (\ref{TwoQuarter1}) with $v=0$ in (\ref{TwoQuarter2}). We employ local/point observers and employ similar optimized dynamic localization based upon the current mesh density function based on (\ref{eq:dynloc}).

We vary the number of spatial mesh points $N_u = 1861, 5101, 20201$ corresponding to cartesian grids of size $30\times 30$, $50\times 50$, and $100\times 100$, respectively, and the time between observations $\Delta = 0.1, 1.0.$ For the look ahead mesh strategy we vary 
the size of the small ensemble with $N_{se} = 2, 4$.
We measure using normalized CPU times per unit of time in which each experiment is based upon $10$ forecast/analysis cycles after discarding the first $10$ cycles. The results obtained are summarized in Table \ref{TableCPU}. Similar RMSE results consistent with the observation error covariance were obtained with the method based upon dynamically updating the individual ensemble meshes based on \cite{KHMVVZ22} (Method A in (\ref{eq:DApred}), (\ref{eq:DAanalysis})) and the method based upon the look ahead strategy (Method B in (\ref{eq:DApred}), (\ref{eq:DAanalysis})) using (\ref{LAHMesh2}) to accumulate the look ahead mesh density function. The CPU times in Table \ref{TableCPU} show that the amount of work required to compute the ensemble forecasts is similar for the two methods. For $N_u=5101$ and $\Delta t=1$, the two methods have similar computational cost while increasing $N_u$ to $20201$ we find that the new method, Method B, is computationally less expensive especially with small ensemble size $N_{se}=2.$ We note that for $\Delta t=0.1$ Method A is uniformly less expensive with the relative difference becoming smaller as $N_u$ increases. We also note that for $N_u=20201$ and $\Delta t = 0.1$ that for Method B the computational cost with $N_{se}=2$ is approximately $10\%$ greater than with $N_{se}=4.$ While a smaller value for the small ensemble size, $N_{se},$ decreases the expense in forming the look ahead mesh, increasing $N_{se}$ has the potential for improved time stepping due to a better conditioned look ahead mesh. 

\begin{table}[]
\begin{center}
\begin{tabular}{||c | c | c | c ||} 
 \hline
 $(N_u,\Delta t)$ & Method B: $N_{se} = 2$ & Method B: $N_{se} = 4$ & Method A \\ [0.5ex] 
 \hline\hline
 $(1861, 1.0)$ & 2.5 & 2.5 & 1 \\ 
 \hline 
 $(5101, 1.0)$ & 5.9 & 6.4 & 5.6  \\ 
 \hline
 $(20201, 1.0)$ & 23.9 & 30.2 & 36.9 \\
 \hline\hline
 $(1861, 0.1)$     & 14.2 & 15.3 & 3.4  \\
 \hline
 $(5101, 0.1)$ & 34.6 & 51.0 & 18.4  \\
 \hline
 $(20201, 0.1)$ & 153	& 134 & 100  \\ 
 \hline
\end{tabular}
\end{center}
    \caption{Normalized CPU times per unit time for ensemble forecasts only comparing Method A and Method B employing small ensemble sizes of $N_{se}=2,4.$}
    \label{TableCPU}
\end{table}

\begin{figure}
    \centering
    \includegraphics[width=6.0true in, trim= 40 270 40 270,clip]{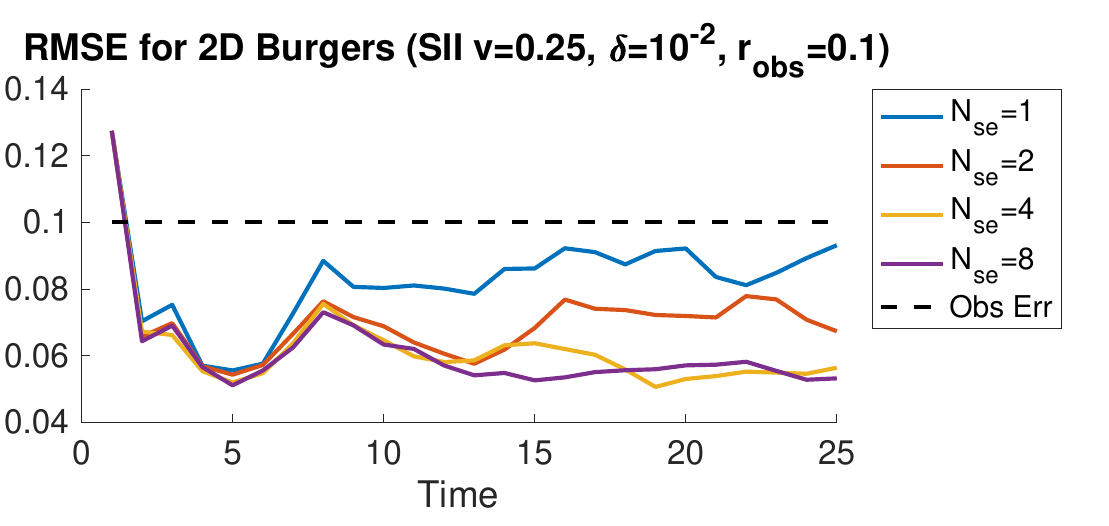}\\
 \includegraphics[width=3.0true in, trim= 160 270 160 270,clip]{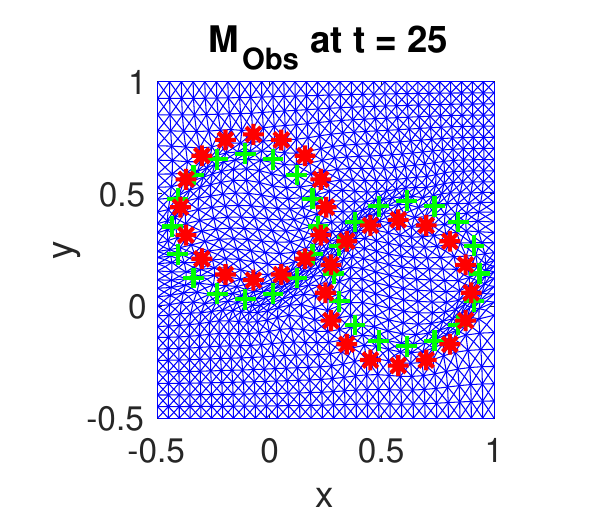}
 \includegraphics[width=3.0true in, trim= 160 270 160 270,clip]{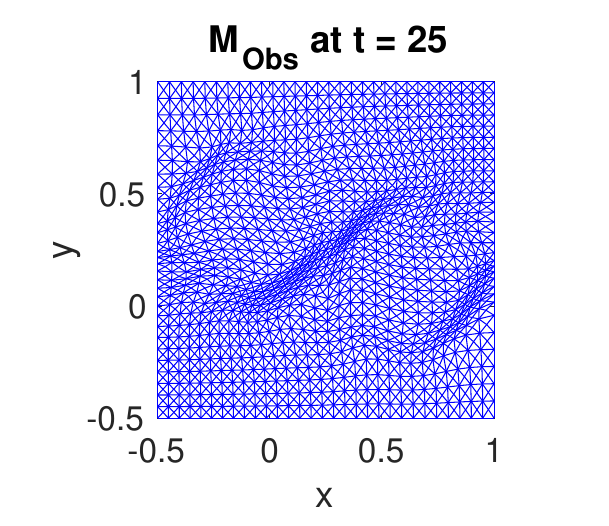}
 \includegraphics[width=4.2true in, trim= 160 280 40 270,clip]{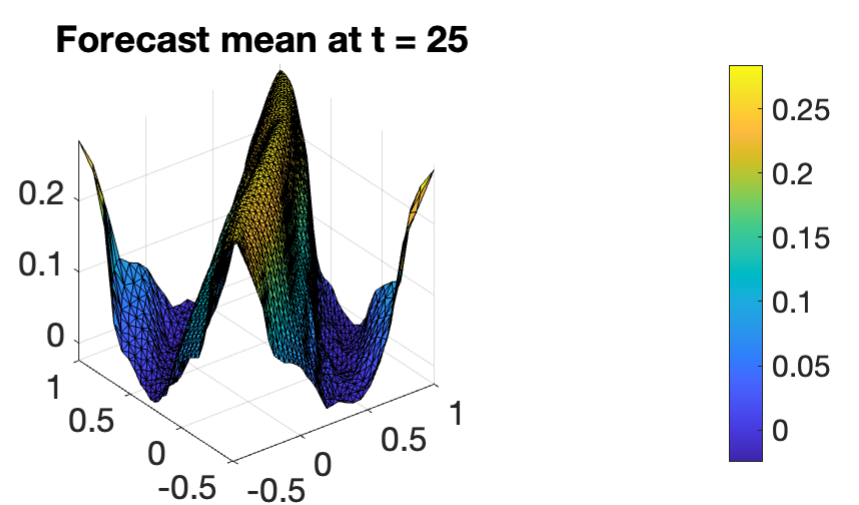}
    \caption{2D Burgers: $N_u=1861$, Remesh after Analysis $\Delta t = 1$}
    \label{fig:Burgers1}
\end{figure}

\begin{figure}
    \centering
    \includegraphics[width=4.0true in, trim= 40 160 40 160,clip]{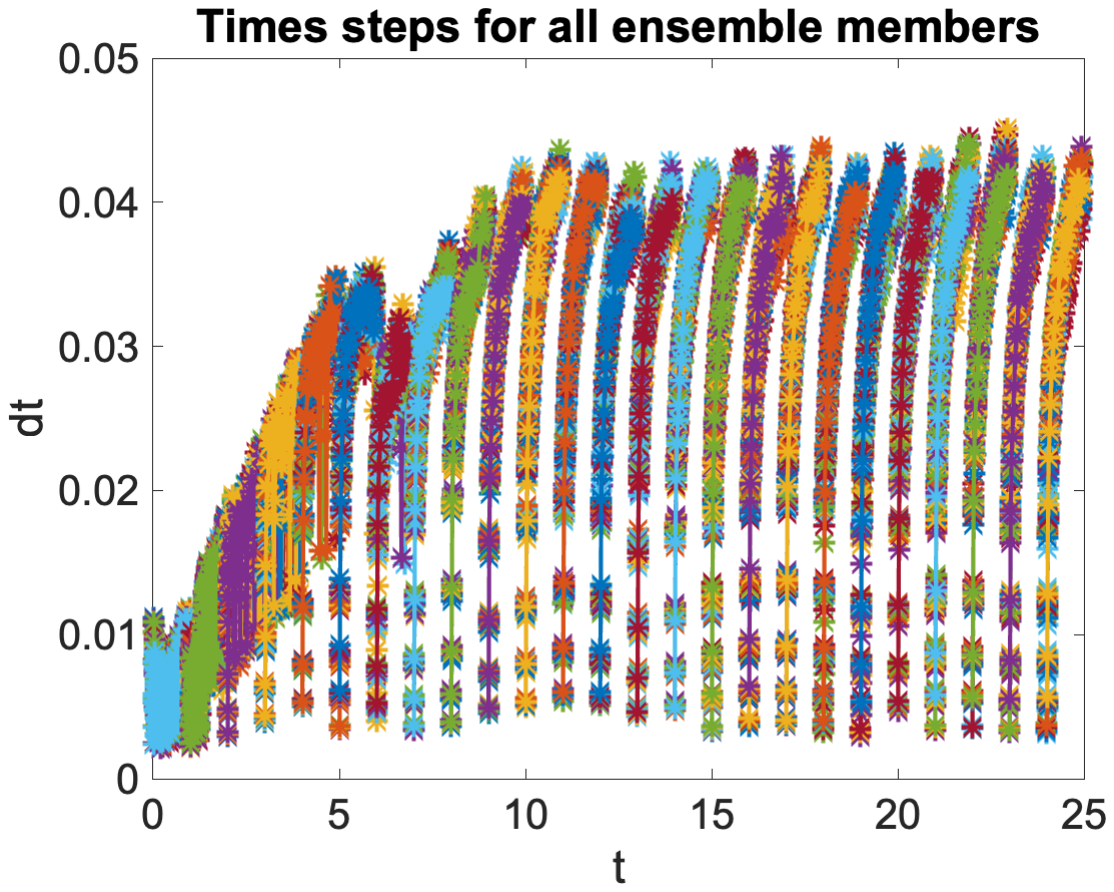}\\
 \includegraphics[width=3.0true in, trim= 60 270 60 270,clip]{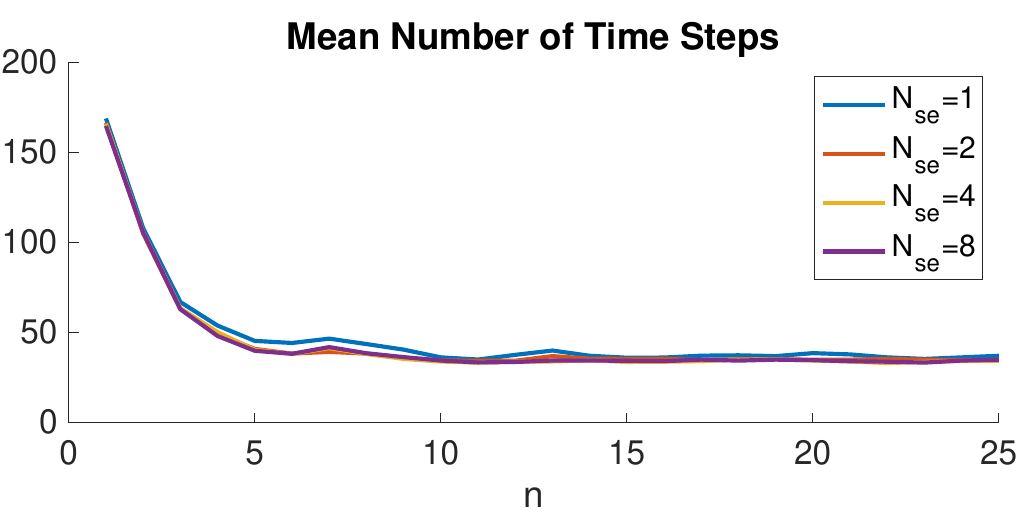}
 \includegraphics[width=3.0true in, trim= 60 270 60 270,clip]{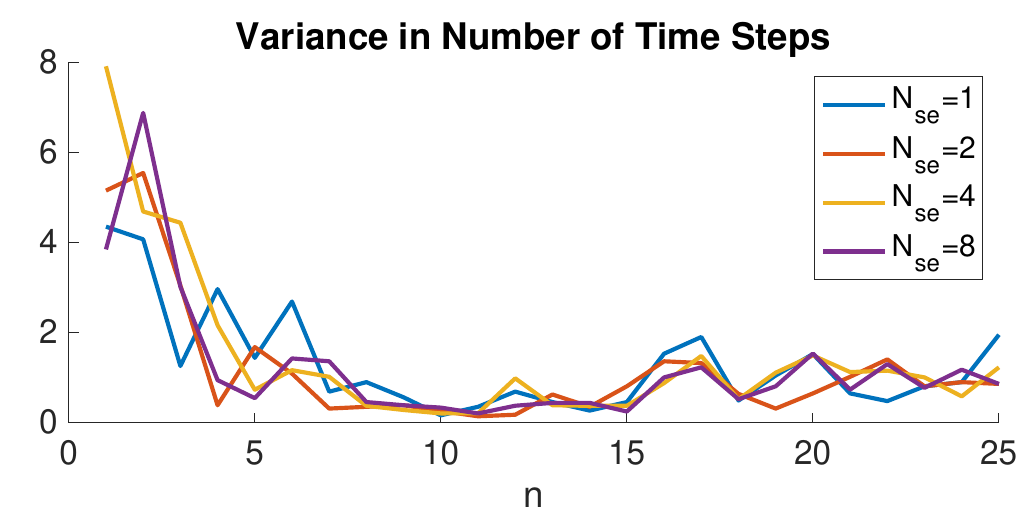}
    \caption{2D Burgers: $N_u=1861$, Remesh after Analysis $\Delta t = 1$}
    \label{fig:Burgers2}
\end{figure}

\begin{figure}
    \centering
    \includegraphics[width=6.0true in, trim= 40 270 40 270,clip]{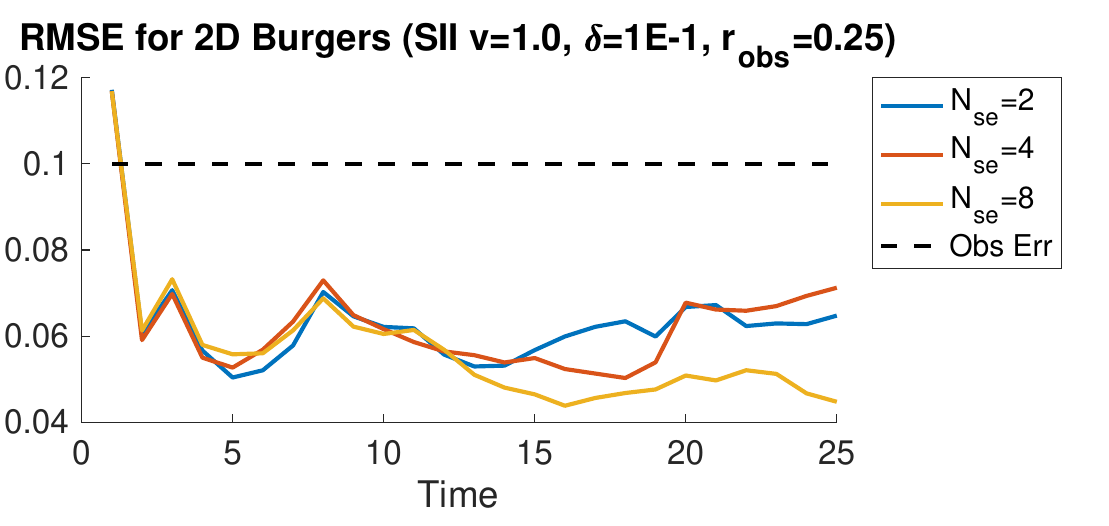}\\
    \caption{2D Burgers: $N_u=1861$, Remesh after Analysis $\Delta t = 1$}
    \label{fig:Burgers3}
\end{figure}

In our next set of experiments we employ covariances
$\bQ=\bR = 10^{-2}\cdot I$ and the observation time scale $\Delta t = 1$ using $N_y=32$ observers split between the two circles in observation Scenario II and uniformly distributed between the two circles. In Figure \ref{fig:Burgers1} we fix the velocity $v=1/4$ in (\ref{TwoCircs}) and vary the size of the small ensemble $N_{se}=1,2,4,8.$ The localization parameter is set to $r_0=1.0$ in (\ref{eq:dynloc}). The $N_y=32$ nonlocal observers are defined using $\delta=10^{-1}$ with radius $B_{r_i}(\hat \bx_i)$ defined with $r_i = r_{obs}=10^{-1}$ in (\ref{HuApproxLoc}). In Figure \ref{fig:Burgers1} we find similar RMSE for $N_{se}=4,8$ while the results for $N_{se}=1,2$ are not as good. We also illustrate the mesh density function $M_{Obs}$ produces a mesh that depends on the form of the ensemble solutions. The results in Figure \ref{fig:Burgers2} show the behavior of the time stepping for the ensemble member forecasts when employing the look ahead mesh. We plot the time steps for all ensemble members between observation times as well as the mean number of time steps taken and variance in the number of time steps between ensemble members. In Figure \ref{fig:Burgers3} we plot the RMSE as a function of time while varying the size of the small ensemble used to form the look ahead mesh. This is for observation scenario II with velocity $v=1$, $N_y=32$ nonlocal observations uniformly distributed over the two circles with $\delta=10^{-1}$ and $r_i=r_{obs}=0.25$ defining the nonlocal observations in (\ref{HuApproxLoc}).
The results are best for $N_{se}=8,$ while for $N_{se}=1$ (using the ensemble mean as the initial condition for the look ahead meshes) the run failed due to small time step.

\begin{figure}
    \centering
    \includegraphics[width=6.0true in, trim= 2 270 48 270,clip]{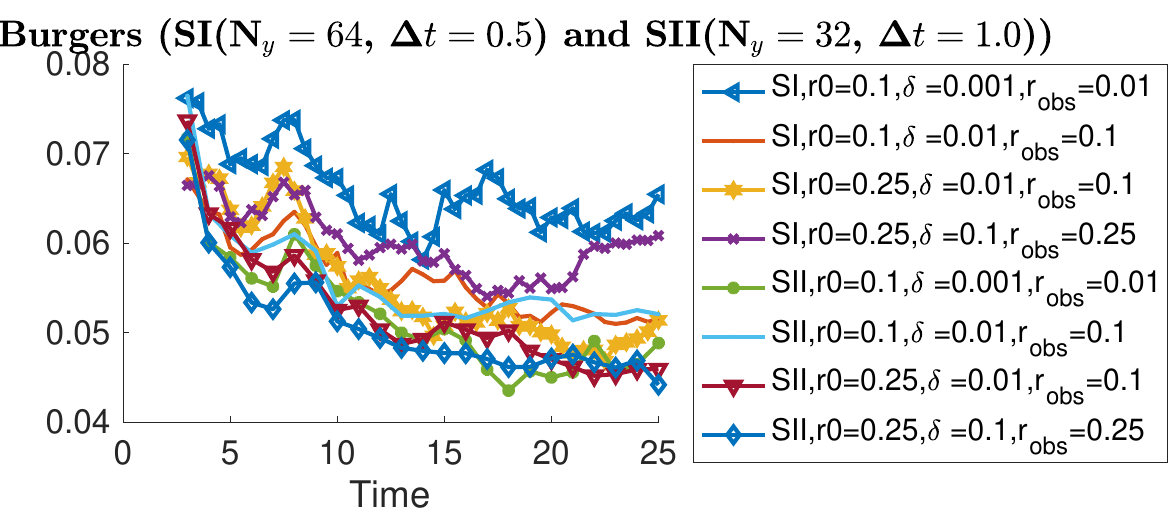}\\
    \caption{2D Burgers: $N_u=1861$, Varying nonlocal observations, localization radius $r_0$.}
    \label{fig:Burgers4}
\end{figure}

We next compare the two different observation scenarios using different observation time scales, nonlocal observation parameters $\delta, r_{obs}$, and localization parameter $r_0.$ We now employ covariances
$\bP_0^b=\bQ= 10^{-1}\cdot I$ and $\bR = 10^{-2}\cdot I$, set the size of the full and small ensemble to $N_e=20$ and $N_{se}=8$, and set the velocity of the moving, nonlocal observations to $v=0.25$. In Figure \ref{fig:Burgers4}  we fix the number of observations $N_y$ and the observation time scale $\Delta t$ so the two observation scenarios have similar skill and compare using nonlocal observers. Stable results are obtained for all values tested. The results for Scenario II with a larger $\Delta t$ and fewer observations $N_y$ are consistently better than those obtained with smaller $\Delta t$ and larger $N_y$ with Scenario I.

%

%

\subsection{Example 2: 2D Shallow Water Equation (Gravity Waves)}

%

Our second example is a SWE on a rectangular domain, configured to approximate gravity waves.
The SWE is frequently used in science and engineering applications to model free-surface flows where the depth is small compared to the horizontal scale(s) of the domain.
The governing equations for this system are
\begin{equation}
\label{eq:swe}\begin{aligned}
  \frac{\partial u}{\partial t} &= \left(-\frac{\partial u}{\partial y}+f \right)v-\frac{\partial}{\partial x}\left(\frac{1}{2}u^2+gh \right)+\nu\nabla^2 u-c_b u, \\
  \frac{\partial v}{\partial t} &= -\left(\frac{\partial v}{\partial x}+f \right)u-\frac{\partial}{\partial y}\left(\frac{1}{2}v^2+gh \right)+\nu \nabla^2 v-c_b v, \\
  \frac{\partial h}{\partial t} &= -\frac{\partial}{\partial x}((h+\underline{h})u)-\frac{\partial}{\partial y}((h+\underline{h})v).
\end{aligned}
\end{equation}
The domain is periodic in the $x$ direction and we employ Dirichlet boundary conditions on the upper and lower boundaries using $v=0$ and $h$ fixed at its initial values. 
Here  $u$ and $v$ are the \(x\)- and \(y\)- components of the   velocity field. The total height of the water column is $h+\underline{h}$, where $h$ is the height of the wave, and $\underline{h}$ is the depth of the ocean, although we employ the flat orography $\underline{h}\equiv 0$. 
 The parameter $g$ is the gravitational constant, $f$ is the Coriolis parameter,
 $c_b$ the bottom friction coefficient, and $\nu$ is the viscosity coefficient.

\begin{figure}
    \centering
    \includegraphics[width=1.45true in, trim= 240 270 240 270,clip]{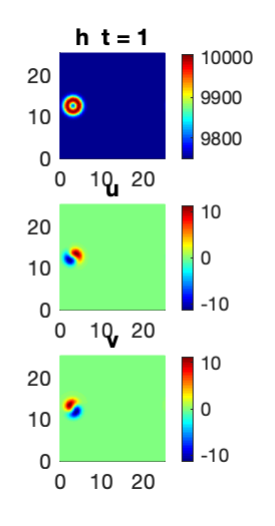}
    \includegraphics[width=1.45true in, trim= 240 270 240 270,clip]{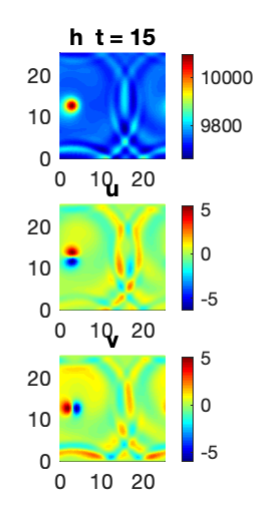}
    \includegraphics[width=1.45true in, trim= 240 270 240 270,clip]{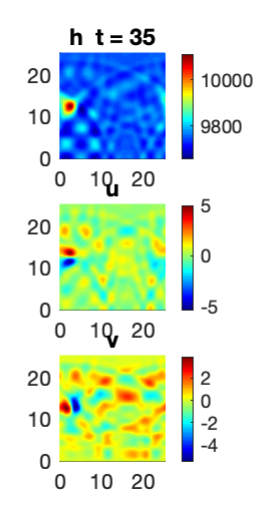}
    \includegraphics[width=1.45true in, trim= 240 270 240 270,clip]{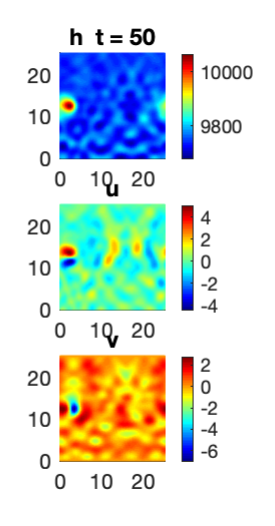}
    \caption{Snapshots of the time evolution of the approximate gravity waves}
    \label{fig:SWE1}
\end{figure}

To solve the SWE, we first non-dimensionalize by scaling space and time, and transform the dependent variables ($h,u,v\to \tilde h, \tilde u, \tilde v$). The physical spatial domain is of the form $(0,L)\times (0,L)$ with $L=25,200$ (km) and we transform to the computational domain $(0,1)\times (0,1)$ and rescale time as $t\to t/T$ with $T=252.$ The dependent variables are transformed as $\tilde h = h/10^4$, $\tilde u = u\cdot h/10^6,$ and
$\tilde v = v\cdot h/10^6.$ In Figure \ref{fig:SWE1} we illustrate the behavior of the true solution we target with snapshots at times $t=1, 15, 35, 50.$

We employ a single mesh to support all three components based upon the mesh density function
\begin{equation}\label{SWEMDF}
\M = \M_1 \cap \M_2,\,\,\, \M_1 = \M(\tilde h),\,\,\, \M_2 = \M(g\cdot \tilde h+ \frac{1}{2}(\tilde u^2 + \tilde v^2)/\tilde h^2).
\end{equation}
We consider a simple observation scenario where all variables $(h,u,v)$ are observed on a uniform grid, see \cite{ABCG-KIMMVV22,KHAN2022,PJBC2022} for other observation scenarios employed for different DA techniques.
In our experiments we employ covariances
$\bQ = \bR =  10^{-2}\cdot I,$ an ensemble of size $N_e=10, 20,$ and observation time scale  of $\Delta t = 0.5$ (hour).

\begin{figure}
    \centering
    \includegraphics[width=6.0true in, trim= 20 270 40 270,clip]{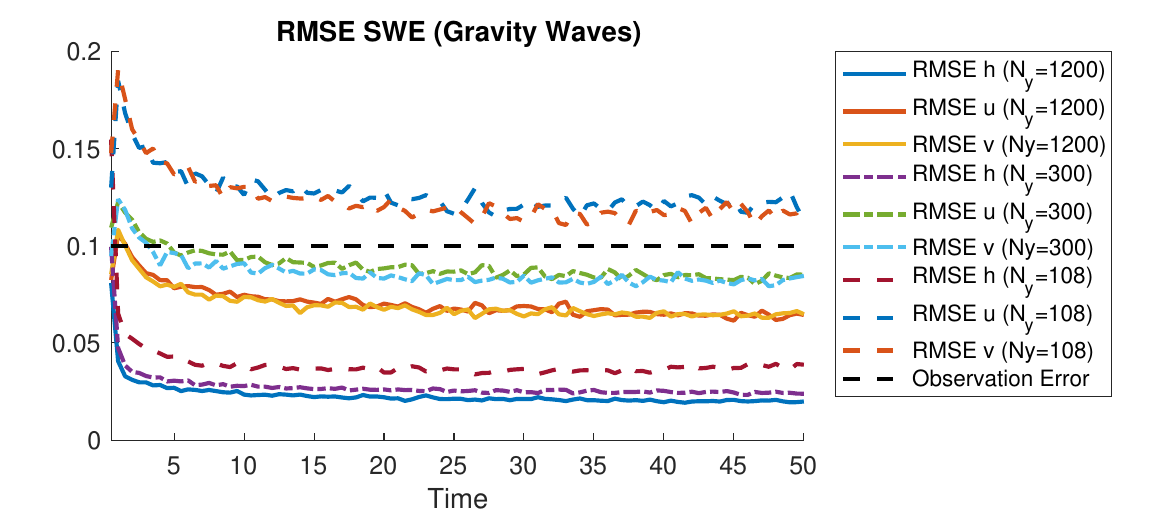} 
    \caption{2D SWE: $N_u=3\cdot 5101$, Point Observations, Observation time scale $\Delta t = 1/2$; Small ensemble size $N_{se}=4$.}
    \label{fig:SWE2}
\end{figure}

In Figure \ref{fig:SWE2} we consider pointwise observations of all components on a uniform grid with $N_y = 3\cdot 20^2, 3\cdot 10^2, 3\cdot 6^2$. The prior covariance and model and data covariances are set to 
$P_0^b = 0.5\cdot I$, $\Sc,\Rc = 10^{-2}\cdot I$. We employ $N_e=10$ ensemble members, form the look ahead mesh using a small ensemble of $N_{se}=4$ and set the observation time scale $\Delta t = 0.5.$ We observe stable results and improved RMSEs as the number of observations is increased.

\begin{figure}
    \centering
    \includegraphics[width=6.0true in, trim= 20 270 40 270,clip]{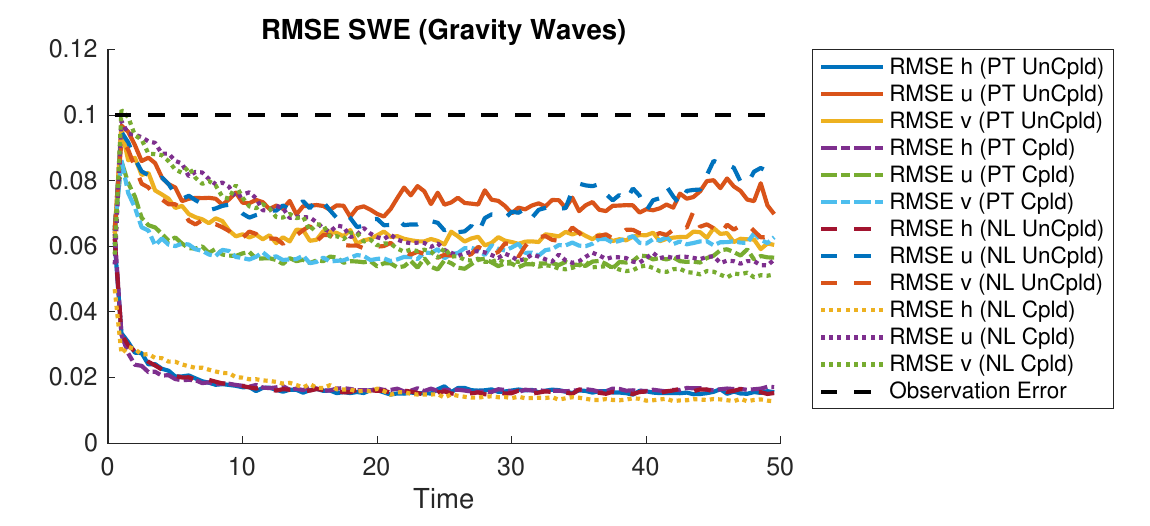} 
    \caption{2D SWE: $N_u=3\cdot 5101$, $N_y=1200$, for nonlocal Observations $(\delta=2.5\cdot 10^{-2}, r_{\tt obs}=5\cdot 10^{-2}),$ Observation time scale $\Delta t = 1/2$; Small ensemble size $N_{se}=8$.}
    \label{fig:SWE3}
\end{figure}

In Figure \ref{fig:SWE3} we consider both local, point observations denoted by {\tt PT} and nonlocal observations denoted by {\tt NL}. We consider both uncoupled and coupled analysis updates.  In the uncoupled case denoted by {\tt UnCpld} the individual variables are updated separately with no forecast covariance between the different variables. In the coupled case denoted by {\tt Cpld} the local analysis for the three dependent variables is performed simultaneously and takes into account forecast covariances between the different variables. What we observe from the experiment is that the height variable $h$ is largely unaffected by the choice of local/nonlocal, coupled/uncoupled, while the velocity variables $u,v$ have improved RMSEs in the coupled case. We also observe the difference of the RMSEs between the variables $u,v$ is smaller in the cases of coupled DA.


\subsection{Example 3: Coupled 1D Kuramoto-Sivashinsky Equations}

We consider a coupled system of two Kuramoto-Sivashinsky equations (KSE) in one space dimension, similar to a system  employed in \cite{coupledKSE}. 
By adjusting the spatial length scale a fast-slow system is created in analogy with the atmosphere-ocean system.
We consider the so-called differentiated form with a Burgers type nonlinearity $uu_x$ and adjust the spatial length scale (width of the spatial domain) and a viscosity like term. The specific form we consider is 
\begin{equation}\label{coupledKSE}
\begin{array}{rl}
   & u_t + uu_x + u_{xx} + \mu_1 u_{xxxx} = c_1 (v-u),\,\,\, 0<x<L_1,\\
   & v_t + vv_x + v_{xx} + \mu_2(x,t) v_{xxxx} = c_2(u-v),\,\,\, 0<x<L_2,
\end{array}
\end{equation}
together with Dirichlet type boundary conditions $u(0,t)=v(0,t)=u(L_1,t)=v(L_2,t)=0$ and
$u_{xx}(0,t)=v_{xx}(0,t)=u_{xx}(L_1,t)=v_{xx}(L_2,t)=0$.

Even in one space dimension the well-posedness of the KSE for non-periodic boundary conditions is not well established (see \cite{larios2014global}). We employ homogeneous Dirichlet boundary conditions and the spatial interval $(0,\pi).$ By considering the odd extension to $(-\pi,0)$, we can identify with odd functions satisfying periodic boundary conditions, similar to a Fourier expansion with a {\tt sine} series. 

We take $L_1 = 3\pi/2$ and $L_2 = 4\pi$ and rescale the spatial variable $x$ so that for both equations $0<x<1$. We fix the viscosity in the first equation as $\mu_1 = 2.5\cdot 10^{-3}.$
The term $\mu_2(x,t)$ depends on space and time and the form we employ is motivated by the ``storm tracker'' versions of the Lorenz '96 equation employed in \cite{StormTrackL96,JMEVV2024a} 
\begin{equation}
    \mu_2(x,t) = \mu_{\tt min} + (\mu_{\tt max} - \mu_{\tt max})(1-|\sin(z)|)^2 ,
\end{equation}
where we take $\mu_{\tt min}= 2.5\cdot 10^{-3}$,  $\mu_{\tt max}= 5\cdot 10^{-2}$, and
$z = \pi (x + \omega\sin(2\pi t))$ with
$\omega = 0.2$.

We evolve KSE by using nonuniform finite differences in space and a linearly implicit backward Euler method in time. 
    We achieve a linearly implicit method by splitting $(u^2/2)_x = uu_x$ as
    $(uu_x)_{n+1} \approx u_n(u_x)_{n+1} + u_{n+1}(u_x)_n - u_n(u_x)_n$.
This results in the following backward Euler discretization for the $u$ and $v$ components,
    \begin{align}
     u_{n+1} + \Delta t \bigg(&\frac{1}{L_1}\left(u_n(u_x)_{n+1} + u_{n+1}(u_x)_n\right)+\frac{1}{L_1^2}(u_{n+1})_{xx} \nonumber \\ +\frac{\mu_1}{L_1^4}&(u_{n+1})_{xxxx}+c_1(u_{n+1}-v_{n+1})\bigg) = u_{n} +  \frac{\Delta t}{L_1} u_n(u_x)_n,\\
     v_{n+1} + \Delta t \bigg(&\frac{1}{L_2}\left(v_n(v_x)_{n+1} + v_{n+1}(v_x)_n\right)+\frac{1}{L_2^2}(v_{n+1})_{xx} \nonumber \\ +\frac{\mu_2}{L_2^4}&(v_{n+1})_{xxxx}+c_2(v_{n+1}-u_{n+1})\bigg) = v_{n} +  \frac{\Delta t}{L_2} v_n(v_x)_n,
    \end{align}
with a complete discretization obtained by substituting nonuniform finite difference approximations, see Appendix B. The mesh is updated after the analysis by generating a mesh using (\ref{MeshUpdate}) and one of the accumulated metric tensors (\ref{LAHMesh1}), (\ref{LAHMesh2}) to produce the next robust, nonuniform common mesh. The ensemble members are interpolated from the previous common mesh to the new common mesh using linear interpolation. 





\begin{figure}
    \centering
    \includegraphics[height=2.0true in,width=3.0true in, trim= 20 270 20 270,clip]{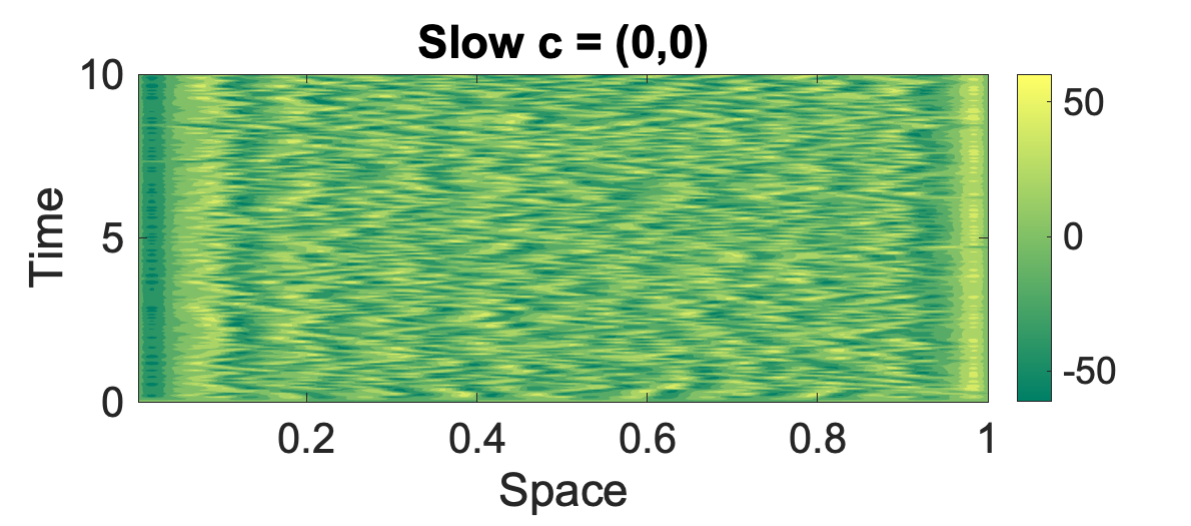}
    \includegraphics[height=2.0true in,width=3.0true in, trim= 20 270 20 270,clip]{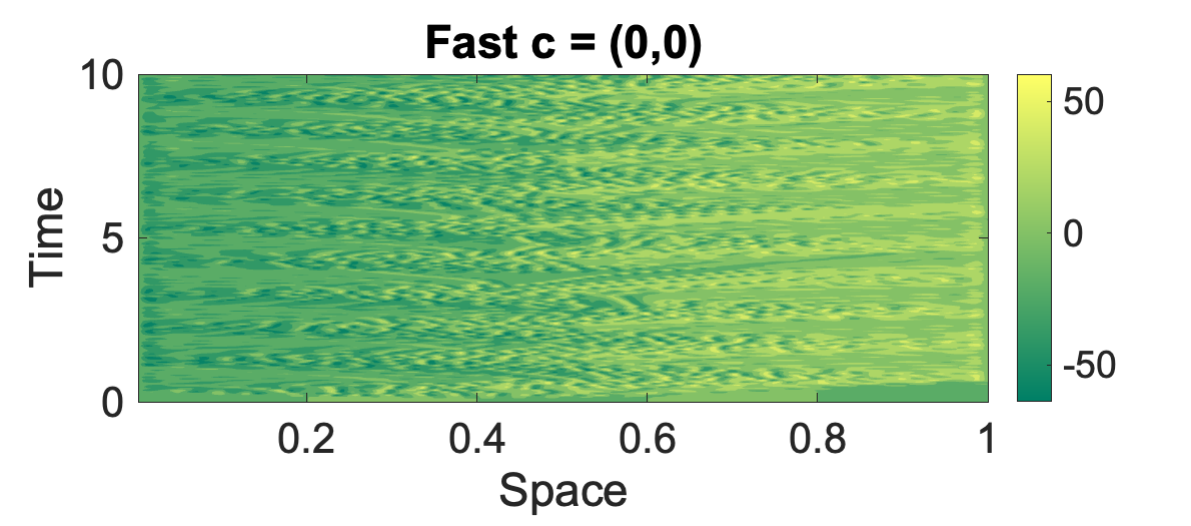}
    \caption{Space time plots of solution to (\ref{coupledKSE}) without coupling $(c_1,c_2)=(0,0).$}
    \label{fig:SpaceTimeKSE}
\end{figure}

Figure \ref{fig:SpaceTimeKSE} illustrates the solution behavior of the uncoupled KSE system over the time interval of $[0,10].$ With increased coupling similar non-equilibrium solution dynamics is observed, while for coupling larger than $c=(1/2,1/2)$ we observed equilibrium or near equilibrium behavior.


\begin{figure}
    \centering
    \includegraphics[width=0.65\linewidth]{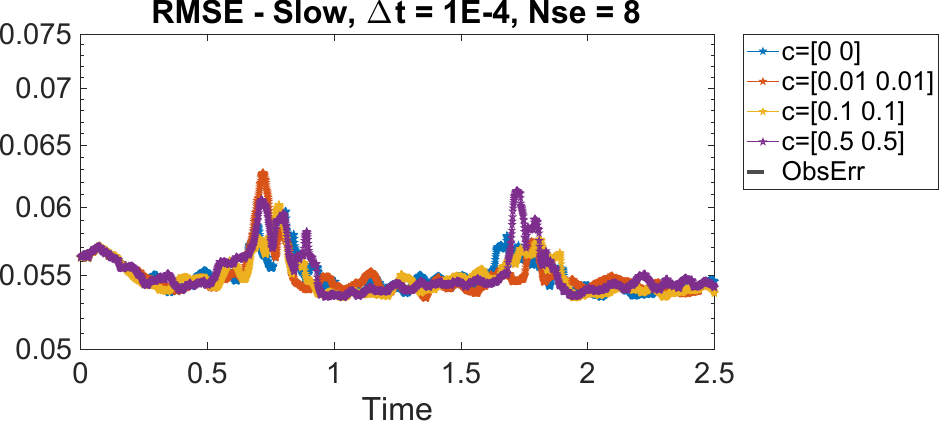}
    \includegraphics[width=0.65\linewidth]{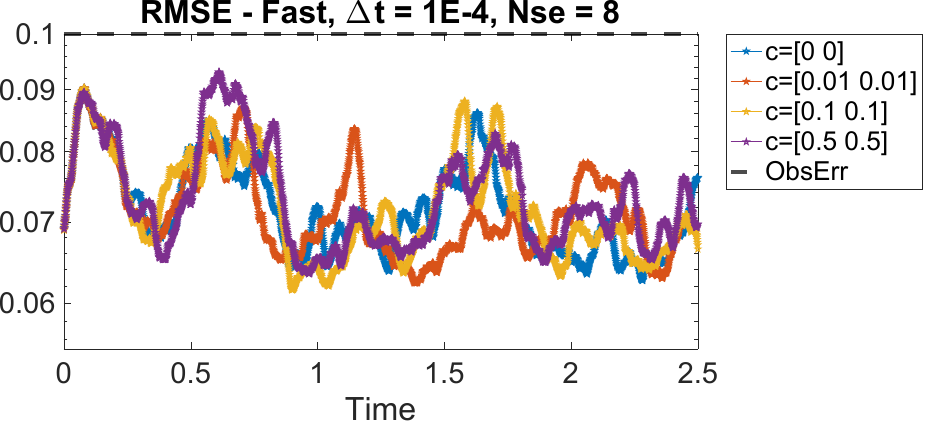}
    \caption{RMSE of KSE with observation time scale $\Delta t = 10^{-4}$, small ensemble size $N_{se}=8$, and PDE coupling values of $[0,0],[0.01,0.01],[0.1,0.1], \text{ and } [0.5,0.5]$. Plotted using moving averages over a time interval of $0.05$.}
    \label{fig:KSE15}
\end{figure}

In Figure \ref{fig:KSE15} we compare different coefficients of the coupling terms between the two KS equations (\ref{coupledKSE}). We use $400$ uniformly space, static observations on the spatial domain of each component. We use coupled DA, which includes the covariances between the components during the analysis update. We use the ``Non-Flow" form of the metric tensor (\ref{LAHMesh2}) for incorporating the observations into the look ahead mesh. We employ a timestep between the observations of $\Delta t = 10^{-4}$ and a computational timestep for the PDE solver of $h=\Delta t/4$. The size of the full and small ensemble are $N_e = 20$ and $N_{se} = 8$. We use covariances $P^b_0=Q=0.1\cdot I$, $R = 0.01\cdot I$. The RMSEs obtained in Figure \ref{fig:KSE15} are relatively uniform across this range of PDE coupling coefficients. 

%
\begin{figure}
    \centering
    \includegraphics[width=0.65\linewidth]{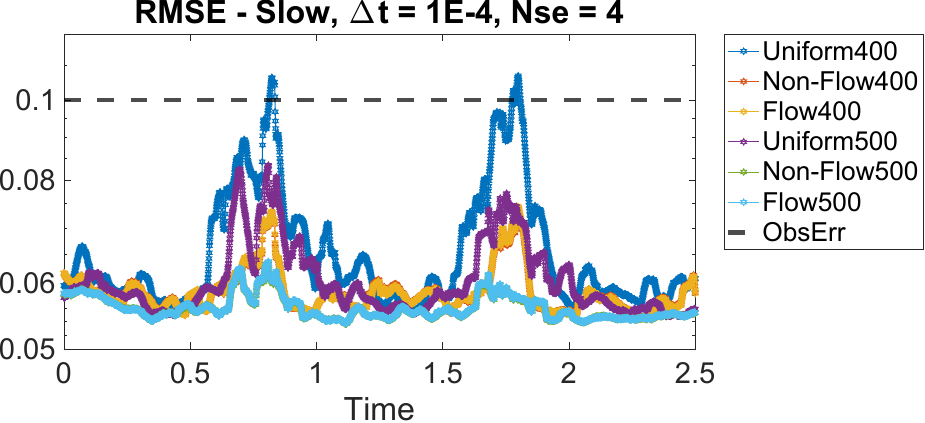}
    \includegraphics[width=0.65\linewidth]{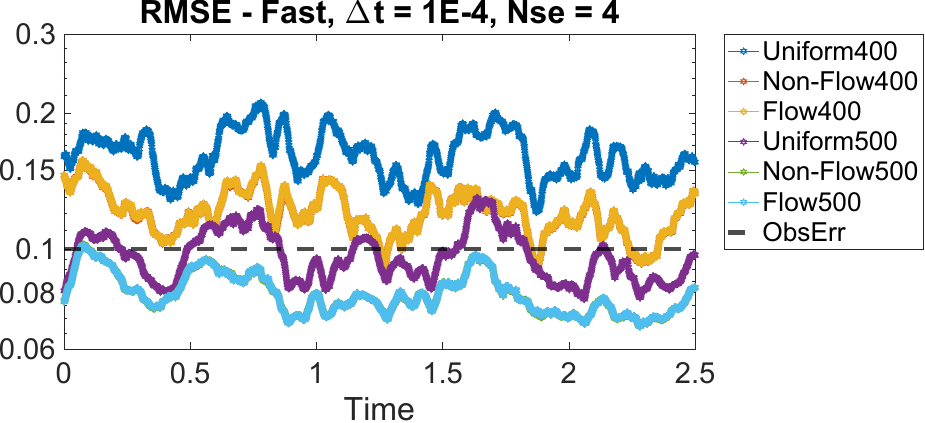}
    \caption{KSE: RMSE of DA on a static uniform mesh of 400 or 500 intervals and an adaptive moving mesh of 401 or 501 points using the flowing and not flowing metric tensor. Observations are located at the midpoints of every other interval.  Plotted using moving averages over a time interval of $0.05$. Control and Flow yield nearly identical values.}
    \label{fig:KSE10}
\end{figure}

In Figure \ref{fig:KSE10} we consider a uniform static mesh of 500 intervals, $\chi_{ref} = \{x_k=k\cdot\Delta x: k=0,...,500\}$ with $\Delta x = 1/500$ and observations placed at $X_{obs} = \{x_k^o = x_0^o + k\cdot \Delta x_{obs}: k = 0,..,249\}$ with $x_0^o =10^{-3}$ and $\Delta x_{obs}=1/250$, and an adaptive mesh of 501 mesh points and the same observations using both metric tensor scenarios. We perform a similar experiment using a uniform static mesh of 400 intervals in which $\chi_{ref}$ is formed using $\Delta x = 1/400$ and the observations are located on $X_{obs}$ using 
$x_0^o=2.5\cdot 10^{-3}$ and $\Delta x_{obs}=1/200$
together with an adaptive, nonuniform mesh of 401 mesh points and the same observation locations.
We test using both the ``Non-Flow'' and ``Flow'' accumulated metric tensors (\ref{LAHMesh2}) and (\ref{LAHMesh1}), respectively.
With the use of static (time independent) observations we observe small difference in the choice of the accumulated metric tensor but the use of an adaptive, nonuniform mesh consistently outperforms the corresponding uniform mesh.

\begin{figure}


        {\centering
    \includegraphics[scale=0.53]{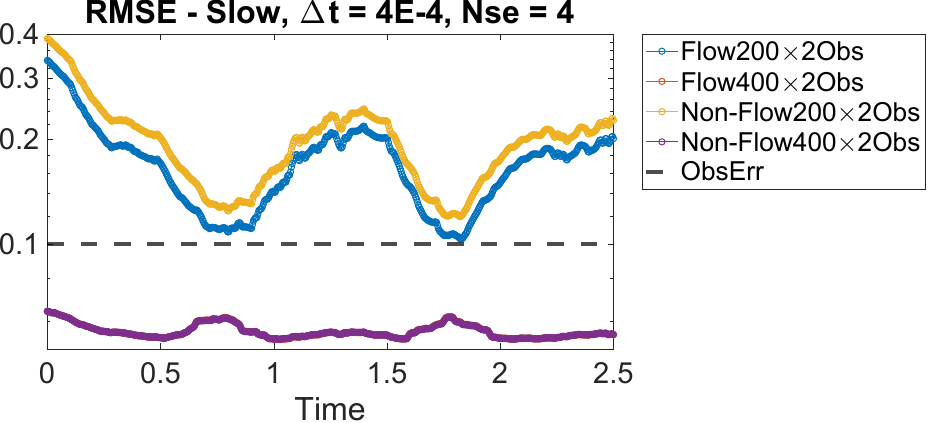}
    \includegraphics[scale=0.53]{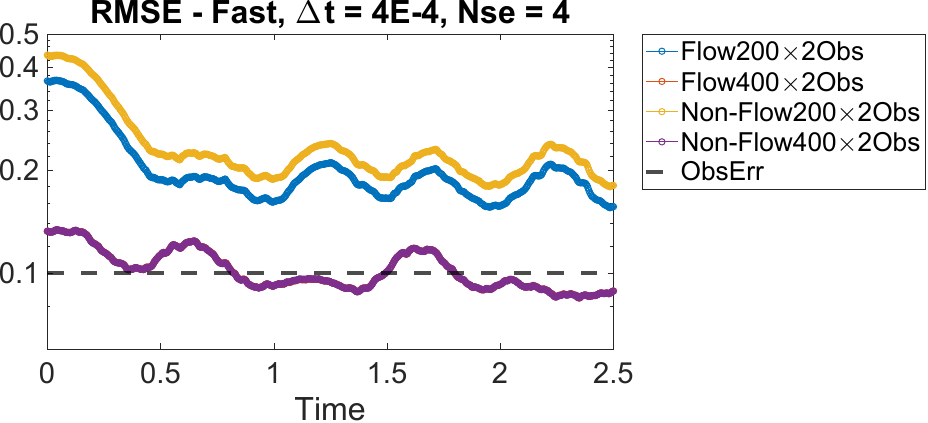} } 
    
    
    
    

    
    \caption{KSE: RMSE of 200$\times$2 moving observations vs 400$\times$2 static observations. Post Look Ahead MT intersection vs Flowing MT intersection. Plotted using moving averages over a time interval of $.05$.
    }
    \label{fig:KSE12}
\end{figure}

\begin{figure}

        {\centering
    \includegraphics[scale=0.55]{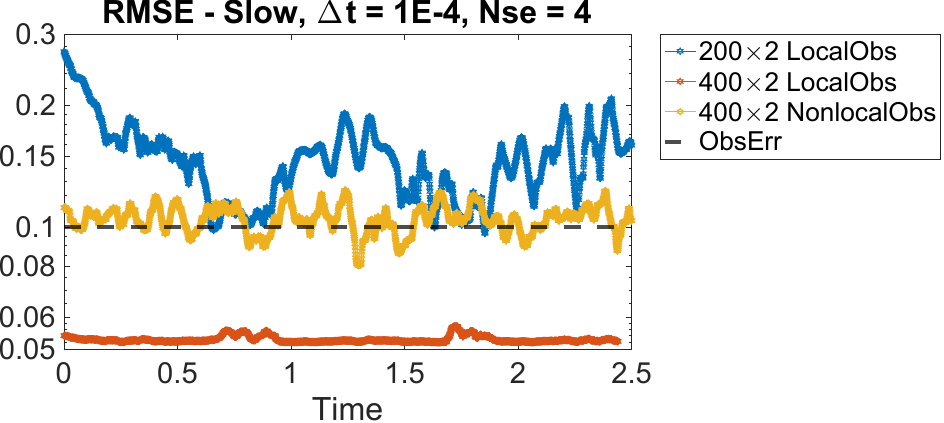}
    \includegraphics[scale=0.55]{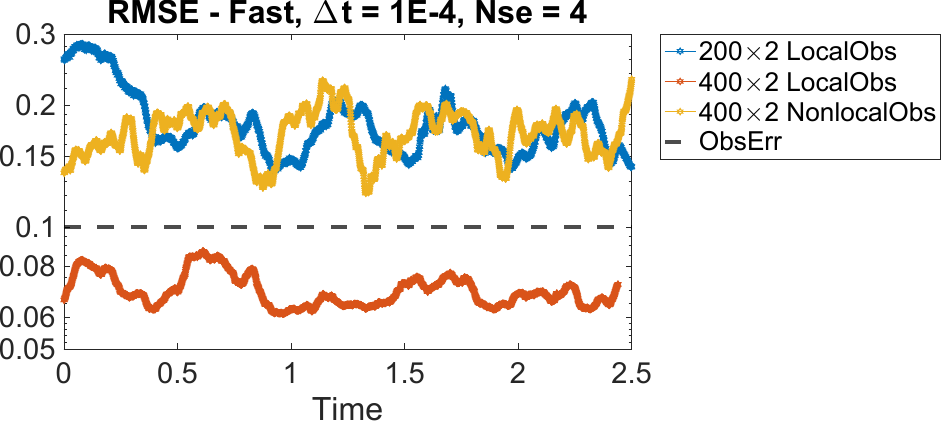} } 
    
    \caption{KSE: RMSE of 200$\times$2 moving local observations, 400$\times$2 static local observations, 
    and 400$\times$2 static nonlocal observations. Plotted using moving averages over a time interval of $.05$.}
    \label{fig:KSE13}
\end{figure}

In Figures \ref{fig:KSE12} and \ref{fig:KSE13} we compare cases with different types of observations,
both static and moving as well as local and nonlocal observations. We also test the use of the ``Non-Flow'' and ``Flow'' accumulated metric tensors (\ref{LAHMesh2}) and (\ref{LAHMesh1}), respectively, with nonlocal observations using (\ref{HuApproxLoc}) with $\delta = 10^{-3}$ and
$r_{\tt obs} = r_i = 2\cdot 10^{-2}$ for all observation locations.
For the case of moving observations, we consider time dependent observation locations
$x_k^o(t) = x_k^o(0) + \frac{1/4}\sin(\omega \pi t)$
for $k=1,...,200\equiv N_y/2$ where we take $\omega=750$ and the initial locations $x_k^o(0)$ are uniformly distributed on the interval $[1/4,3/4].$
With $N_y/2 = 400$ static observations, we consider observation locations uniformly distributed over $[0,1]$, while for static nonlocal observation, we use the same locations uniformly distributed over the spatial domain  with a nonlocal observation operator (\ref{HuApproxLoc}). We use the ``Non-Flow" form of the metric tensor for the observations. We employ a timestep between the observations of $\Delta t = 10^{-4}$ and a timestep for the PDE solver of $h=\Delta t/4$. The size of the full and small ensemble are $N_e = 20$ and $N_{se} = 4$. We use covariances $P^b_0=Q=0.1\cdot I$, $R = 0.01\cdot I$. In Figure \ref{fig:KSE12} we see the advantage with moving observations of employing the ``Flow'' technique (\ref{LAHMesh1}). While there is an improvement with moving observations, with non-moving observations both techniques yield essentially identical RMSEs. In Figure \ref{fig:KSE13} we see similar RMSE values with static local and moving nonlocal observations (albiet with a smaller number of local observers) and improved RMSE values with a larger number of static local observations. With the form of the moving observations, the analysis updates at each observation time only update approximately half of the spatial locations. nonlocal observations are only providing local averages as opposed to precise local values with local observers which translates into improved RMSEs.

\section{Conclusion and Future Work}\label{Conclusion}

In this work we continue the development of a flexible, modular framework for DA with adaptive meshing. Using a mesh density or metric tensor formulation, we develop robust meshes that combine the needs for accurate, efficient ensemble forecasts and the flow of information from observational data. We develop robust meshes that support both the ensemble forecasts and the analysis update based on so-called look ahead meshes. This reduces the need for interpolation during the forecast-analysis cycle and since all ensemble meshes are supported on the same mesh increases the efficieny of ensemble forecasts. We consider local and nonlocal observations with both fixed locations and time dependent observers. Using three model problems, two of which are multi-component models and two of which are in two space dimensions, we illustrate the utility of the techniques developed. 




This work suggested several areas for further investigation. Among these are more investigation in the utility of coupled DA techniques for multi-components systems, the combination of local and nonlocal observational data and how to effectively update the analysis. The further development of metric tensors or mesh density functions in a DA context is another area for in depth investigation. These include mesh density functions to reduce different types of DA errors through meshing based on the sensitivity of the DA as well as the ad hoc need for variable resolution to increase DA precision dynamically in different space-time regions. We are also interested in employing so-called mesh quality measures to adaptively determine the size of the small ensemble used to determine robust look ahead meshes.

%
%
%
%

    \bibliographystyle{abbrv}
    \bibliography{KU_MMDAnew}

\newpage

\begin{appendices}

\section{Mesh Density nonlocal Observations}\label{AppA}

We next derive the goal-oriented metric tensor that minimizes the $L^2$ linear interpolation error on the mesh $\T_h$. Define
\begin{equation}
 E_h = \sum_{i=1}^{N_y} \left|\int_\Omega (u - \pi_h u)G(\hat \bx_i - \bx)d\bx   \right|,
\end{equation}
where $\pi_h$ is the linear interpolation operator on $\T_h.$
Then we have 
\begin{align*}
 E_h \leq & \sum_{i=1}^{N_y} \sum_K \int_K |(u - \pi_h u)|\cdot |G(\hat \bx_i - \bx)|d\bx \\
    \approx & \sum_{i=1}^{N_y} \sum_K  |G(\hat \bx_i - \bx_K)|\int_K |(u - \pi_h u)|d\bx \\
    \leq & \sum_{i=1}^{N_y} \sum_K  |K|^{1/2} \|(u - \pi_h u)\|_{L^2(K)}\cdot |G(\hat \bx_i - \bx_K)| \\
    = &  \sum_K  |K|^{1/2} \|(u - \pi_h u)\|_{L^2(K)}\cdot \sum_{i=1}^{N_y} |G(\hat \bx_i - \bx_K)|,
\end{align*}
where $\bx_K$ is the centroid of the element $K$. Use a known estimate \cite{Huang2003}
\begin{equation}
    \|u-\pi_h u\|_{L^2(K)}\lesssim C|K|^{1/2}{\rm tr}((F'_K)^T|H_K|F'_K),
\end{equation}
where $F'_K$ is the Jacobian of the affine mapping $F_K:\hat K \to K$ from the reference element to the given element, and $|H_K| = \sqrt{H_K^2}$ with $H_K$ an approximate Hessian of $u$ on the element $K.$ If we regularize this with a positive parameter $\alpha_h>0$, then
\begin{equation}\label{AKdef}
E_h \lesssim C \alpha_h \sum_K |K| {\rm tr}((F'_K)^T A_K F'_K),\quad
A_K = I + \frac{1}{\alpha_h} |H_K| \sum_{i=1}^{N_y} |G(\hat \bx_i - \bx_K)| .
\end{equation}

We choose the optimal metric tensor $\M$ such that $E_h$ is minimized at uniform meshes in this metric that satisfy
the alignment and equidistribution conditions
\begin{equation}
    \frac{1}{d} {\rm tr}((F'_K)^T\M_K F'_K) = {\rm det}((F'_K)^T\M_K F'_K)^{1/d},\quad 
    |K|{\rm det}(\M_K)^{1/2} = \frac{1}{N} \sigma_h ,
\end{equation}
where $N$ denotes the number of elements and $\sigma_h = \sum_K |K| {\rm det}(\M_K)^{1/2}$.
Comparing the alignment condition with the right-hand side of $E_h$ suggests the choice
$\M_K = \rho_K A_K$ for some scalar $\rho_K$. Then
\[
\frac{1}{d}{\rm tr}((F'_K)^TA_K F'_K) = {\rm det}((F'_K)^TA_K F'_K)^{1/d} = |K|^{2/d} {\rm det}(A_K)^{1/d} .
\]
Using this, we can rewrite $E_h$ as $E_h \lesssim C \alpha_h \sum_K |K|^{1+2/d} {\rm det}(A_K)^{1/d}$.

An optimal form of $E_h$ (e.g., see \cite{Huang2010}) is
\[
E_h \lesssim C \alpha_h \sum_K |K| {\rm det}(\M_K)^{1/2}\cdot (|K|{\rm det}(\M_K)^{1/2})^{\beta} = C \alpha_h \sigma_h (\frac{\sigma_h}{N})^\beta = C \alpha_h \frac{\sigma_h^{1+\beta}}{N^\beta}.
\]
Equating, we obtain $|K|^{1+\beta} = |K|^{1+d/2}$ and ${\rm det}(\M_K)^{(1+\beta)/2} = {\rm det}(A_K)^{1/d}$.
From these and recalling $\M_K = \rho_K A_K$, we get $\beta = d/2$ and
\[
\rho_K = {\rm det}(A_K)^{-\frac{1}{d+2}}, \quad \M_K = {\rm det}(A_K)^{-\frac{1}{d+2}} A_K,
\]
where $A_K$ is given in (\ref{AKdef}). 
The regularization parameter can be chosen (e.g., $\alpha_h = 1$) or an optimal $\alpha_h$ can be derived.


\section{nonuniform Finite Differences for KS equations}\label{AppB}

    We approximate each derivative by nonuniform finite differences that are included here for completeness
    \begin{align}
    u_x^j & \approx \frac{u^{j}-u^{j-1}}{x^{j}-x^{j-1}},
    \\
    u_{xx}^j & \approx \frac{u_x^{j+1}-u_x^j}{x^{j+1}-x^j}
             \approx \frac{2}{x^{j+1}-x^{j-1}}\left(\frac{u^{j+1}-u^{j}}{x^{j+1}-x^{j}}  - \frac{u^{j}-u^{j-1}}{x^{j}-x^{j-1}}
    \right) ,
    \\
u_{xxxx}^j & \approx \frac{2}{x^{j+1}-x^{j-1}}\left(\frac{u_{xx}^{j+1}-u_{xx}^{j}}{x^{j+1}-x^{j}}
    - \frac{u_{xx}^{j}-u_{xx}^{j-1}}{x^{j}-x^{j-1}}
    \right)
    \notag
    \\
    & \approx  \frac{2}{x^{j+1}-x^{j-1}}\Bigg [
    \frac{2}{(x^{j+2}-x^{j})(x^{j+1}-x^j)}\left(\frac{u^{j+2}-u^{j+1}}{x^{j+2}-x^{j+1}}
    - \frac{u^{j+1}-u^{j}}{x^{j+1}-x^{j}}
    \right)\nonumber\\
    & \quad - 
    \frac{2}{(x^{j+1}-x^{j-1})(x^{j+1}-x^j)}\left(\frac{u^{j+1}-u^{j}}{x^{j+1}-x^{j}}
    - \frac{u^{j}-u^{j-1}}{x^{j}-x^{j-1}}
    \right)\nonumber\\
    & \quad - \frac{2}{(x^{j+1}-x^{j-1})(x^{j}-x^{j-1})}\left(\frac{u^{j+1}-u^{j}}{x^{j+1}-x^{j}}
    - \frac{u^{j}-u^{j-1}}{x^{j}-x^{j-1}}
    \right)\nonumber\\
    & \quad + 
    \frac{2}{(x^{j}-x^{j-2})(x^{j}-x^{j-1})}\left(\frac{u^{j}-u^{j-1}}{x^{j}-x^{j-1}}
    - \frac{u^{j-1}-u^{j-2}}{x^{j-1}-x^{j-2}}
    \right)
    \Bigg ] .
\end{align}

\end{appendices}

\end{document}